\documentclass[a4paper,11pt,oneside,notitlepage, draft]{amsart}
\usepackage[dvips]{graphicx}
\usepackage{latexsym}
\usepackage{amsfonts}
\usepackage{amssymb}
\usepackage{amsthm}
\usepackage{amsmath}

\newtheorem{teo}{Theorem}
\newtheorem*{esmpi}{Examples}
\newtheorem*{os}{Remark}
\newtheorem{prop}{Proposition}
\newtheorem{defin}{Definition}
\newtheorem*{notaz}{Notazione}
\newtheorem*{ipts}{Ipotesi}
\newtheorem{lem}{Lemma}
\newtheorem{cor}{Corollary}

    {\begin{teo} \begin{itshape} \hspace{0.1em}}%
    {\end{itshape} \end{teo}}

    {\begin{esmpi} \begin{upshape} \hspace{0.1em}}%
    {\end{upshape} \end{esmpi}}

\newenvironment{definizione}%
    {\begin{defin} \begin{itshape} \hspace{0.1em}}%
    {\end{itshape} \end{defin}}

\newenvironment{osservazione}%
    {\begin{os} \begin{upshape} \hspace{0.1em}}
    {\end{upshape} \end{os}}

\newenvironment{proposizione}%
    {\begin{prop} \begin{itshape} \hspace{0.1em}}
    {\end{itshape} \end{prop}}

    {\begin{notaz} \begin{upshape} \hspace{0.1em}}%
    {\end{upshape} \end{notaz}}

    {\begin{ipts} \begin{upshape} \hspace{0.1em}}%
    {\end{upshape} \end{ipts}}

\newenvironment{lemma}%
    {\begin{lem} \begin{itshape} \hspace{0.1em}}
    {\end{itshape} \end{lem}}

    {\begin{cor} \begin{itshape} \hspace{0.1em}}%
    {\end{itshape} \end{cor}}

\usepackage{syntonly}

\title[]{Nonlinear Schr\"odinger equation on  four-dimensional compact manifolds}

\author{Patrick G\'erard}
\address{Patrick G\'erard:
Unversit\'e Paris Sud,
Math\'ematiques, B\^{a}t 425, 91405 Orsay Cedex, France}
\email{Patrick.gerard@math.u-psud.fr}

\author{Vittoria Pierfelice}
\address{Vittoria Pierfelice: Universit\`a di Pisa, Dipartimento di Matematica,
Via Buonarroti 2, I-56127 Pisa, Italy}
\email{pierfelice@dm.unipi.it}

\begin{document}

\begin{abstract}
We prove two new results about the Cauchy problem in the energy
space for nonlinear Schr\"odinger equations on four-dimensional
compact manifolds. The first one concerns global wellposedness for
Hartree-type nonlinearities and includes approximations of cubic
NLS on the sphere. The second one provides, in the case of zonal
data on the sphere, local wellposedness for quadratic
nonlinearities as well as global wellposedness for small energy
data in the Hamiltonian case. Both results are based on new
multilinear Strichartz-type estimates for the Schr\"odinger group.

 \end{abstract}

\maketitle
%\thispagestyle{empty}
 %\begin{abstract}

 %\end{abstract}

\section{Introduction}

In a recent series of papers ( \cite{BGT4}, \cite{BGT8},
\cite{BGT7},
 see also \cite{BGT3}, \cite{BGTe}) ,
Burq-G\'erard-Tzvetkov investigated the Cauchy problem for
nonlinear Schr\"odinger equations (NLS) on Riemannian compact
manifolds, generalizing the work of Bourgain on tori (\cite{B1},
\cite{B2}). In \cite{BGT4},  Strichartz estimates with fractional
loss of derivatives were established for the Schr\"odinger group.
They led to global wellposedness of NLS on surfaces with any
defocusing polynomial nonlinearity. On three-manifolds, these
estimates also provided global existence and uniqueness for  cubic
defocusing NLS, but they failed to prove the Lipschitz continuity
of the flow map on the energy space. These results were improved
in
 \cite{BGT8}, \cite{BGT7} for specific manifolds such as spheres,
taking advantage of new multilinear Strichartz inequalities for
the Schr\"odinger group (see also \cite{BGT6}). In particular, on
such three-manifolds the Lipschitz continuity and the smoothness
of the flow map on the energy space were established for cubic
NLS, as well as global existence on the energy space for every
defocusing subquintic NLS. \vskip0.5cm \noindent However, none of
the above methods provided global wellposedness results in the
energy space for NLS on four-dimensional manifolds. This is in
strong contrast with the Euclidean case (see \cite{GV2}, \cite{K},
\cite{C}, \cite{TV}). The only available global existence  result
on a compact four-manifold seems to be the one of Bourgain in
\cite{B2}, which concerns defocusing nonlinearities of the type
$|u|u$ and Cauchy data in $H^s({\mathbb T}^4)$, $s>1$. Let us
discuss briefly the reasons of this difficulty. On the one hand,
Strichartz estimates of \cite{BGT4} involve a too large loss of
derivative in four space dimension ; typically, for cubic NLS,
they lead to local wellposedness in $H^s$ for $s>3/2$, which is
not sufficient in view of the energy and $L^2$ conservation laws.
Moreover, these estimates are restricted to $L^p_tL^q_x$ norms
with $p\geq 2$ and the admissibility condition
$$\frac{1}{p}+\frac{2}{q}=1\ ,$$
so that the analysis does not improve when the nonlinearity
becomes subcubic. On the other hand, the analysis based on
bilinear Strichartz estimates is currently restricted to
nonlinearities of cubic type, and on $S^4$ it only yields local
wellposedness in $H^s$ for $s>1$. In fact, this obstruction can be
made more precise by combining two results from \cite{BGT4} and
\cite{BGT8}. Indeed, from Theorem 4 in \cite{BGT4}, we know that
the estimate
$$\int _0^{2\pi }\! \int _{S^4}| e^{it\Delta }f\, (x)|^4\, dt\, dx
\lesssim \, \| f\| _{H^{1/2}(S^4)}^4 $$ is wrong, which, by Remark
2.12 in \cite{BGT8}, implies that the flow map of cubic NLS cannot
be $C^3$ near the Cauchy data $u_0=0$ in $H^1(S^4)$. Moreover,
notice that this phenomenon occurs for zonal data, namely
functions depending only on the distance to a fixed point.
\vskip0.5cm \noindent
 The purpose of this paper is to provide further results on
 four-dimensional manifolds. We shall study two types of NLS
 equations. In section 2, we study NLS with the following nonlocal
 nonlinearity,
\begin{equation}\label{eq1'}
\begin{cases}
&{i} \partial_t u + \Delta u = \left( (1- \Delta)^{-\alpha}  |u|^2 \right) u,\\
&u(0,x)=u_0(x)
\end{cases}
\end{equation}
where $\alpha >0$. Notice that the homogeneous version of this
nonlinearity on the Euclidean space ${\mathbb R}^d$ reads
$$\left (\frac{1}{|x|^{d-2\alpha }}*|u|^2\right )u $$
so that (\ref{eq1'}) can be seen as a variant of Hartree's
equation on a compact manifold.  We obtain the following result.

\begin{teo}\label{t1'}
Let $(M,g)$ be a compact Riemannian manifold of dimension 4 and
let $\alpha >\frac{1}{2}$.  There exists a subspace $X$ of
$\mathcal{C}(\mathbb{R}, H^1(M))$ such that, for every $u_0 \in
H^1(M)$, the Cauchy problem (\ref{eq1'}) has a unique global
solution $u \in X$. Moreover, in the special case $M$ is the
four-dimensional standard sphere $M=\mathbb{S}^4$, the same result
holds for all values $\alpha>0$ of the parameter.
\end{teo}

The proof of Theorem \ref{t1'} relies on the combination of
conservation laws for equation (\ref{eq1'}) with the following
quadrilinear estimates,
\begin{equation*}
\begin{aligned}
\sup _{\tau \in {{\mathbb R}}}&
\left| \int_{{\mathbb R}} \int_{M}
\chi (t)\,  e^{it \tau}
(1-\Delta)^{-{\alpha}}
( u_1 \overline{u}_2){u_3}  \overline{u}_4dx dt  \right|\\
&\leq C (m(N_1,\cdots ,N_4))^{s_0} \|f_1\|_{L^2(M)} \|f_2\|_{L^2(M)}  \|f_3\|_{L^2(M)}  \|f_4\|_{L^2(M)},\\
\end{aligned}
\end{equation*}
 for every $\chi \in {\mathcal C}^\infty
_0({\mathbb R})$, for every $s_0< 1$ and for
$f_1,f_2,f_3, f_4$ satisfying
\begin{equation*}
\mathbf{1}_{\sqrt{1-\Delta} \in [N_j,2N_j]}(f_j) = f_j, \;\;
j=1,2,3,4. \
\end{equation*}
Here and in the sequel $m(N_1,\cdots ,N_4)$ denotes the product of
the smallest two numbers among
 $N_1, N_2,$ $N_3, N_4$. Moreover $u_j$ and $f_j$ are linked by
 $$u_j(t,x)=S(t) f_j(x), \; j=1,2,3,4,$$
 where $S(t)=e^{it\Delta }.$
Notice that, compared to the multilinear estimates used in \cite{BGT7},
 a frequency variable $\tau $ is added to the left hand side of the estimate.
It would be interesting to know if the smallest value of $\alpha $
for which these estimates (and hence Theorem \ref{t1'}) are valid
depends or not on the geometry of $M$.

\vskip0.5cm \noindent

In Section 3, we come back to power nonlinearities. Since we want
to go below the cubic powers and at the same time we want to use
multilinear estimates, we are led to deal with quadratic
nonlinearities. In other words, we study the following equations,

\begin{equation}\label{eq3'}
{i} \partial_t u + \Delta u =  q(u),
\end{equation}
where $q(u)$ is a homogeneous quadratic polynomial in $u,\overline
u$
$$q(u)= a u^2 + b \overline{u}^2+ c|u|^2.$$

We start with wellposedness results.

\begin{teo}\label{t2}
If $(M,g)$ is the four-dimensional standard sphere , then the
Cauchy problem (\ref{eq3'}) is (locally in time) uniformly
well-posed in $H^s_{{\rm zonal}}(S^4)$ for every $s> \frac{1}{2}
$, where $H^s_{{\rm zonal}}(S^4)$ denotes the $H^s$ space of zonal
functions relative to some pole $\omega \in S^4$ : $f(x)=\tilde
f(\langle x, \omega \rangle )\ $.
\end{teo}

 The main
tool in the proof of Theorem \ref{t2}  is the following trilinear
estimate on linear solutions $u_j(t)=S(t)f_j$,
\begin{equation}\label{TE}
\begin{aligned}
\sup _{\tau \in {{\mathbb R}}}&\left| \int_{{\mathbb R}}
\int_{\mathbb{S}^4} \chi (t)\, e^{it \tau} {\mathcal T}(\;
u_1(t,x),\; u_2(t,x),\;
u_3(t,x))\; dx\, dt  \right| \\
& \leq C\,(\min(N_1,N_2, N_3))^{s_0} \|f_1\|_{L^2(\mathbb{S}^4)} \|f_2\|_{L^2(\mathbb{S}^4)}  \|f_3\|_{L^2(\mathbb{S}^4)},\\
\end{aligned}
\end{equation}
for every ${\mathbb R}$-trilinear expression ${\mathcal T}$ on
${\mathbb C}^3$, for every $\chi \in {\mathcal C}^\infty
_0({\mathbb R})$, for every $s_0>1/2$ and for zonal functions
$f_1,f_2,f_3$ satisfying
\begin{equation*}
\mathbf{1}_{\sqrt{1-\Delta} \in [N_j,2N_j]}(f_j) = f_j, \;\;
j=1,2,3 \ .
\end{equation*}
It would be interesting to know whether the above
estimate holds
 with non zonal functions for some $s_0<1$ ; this would extend the above theorem to
 any finite energy Cauchy data.

 \noindent Notice that a subclass of these equations consists of Hamiltonian
 equations
$$q(u)=\frac{\partial V}{\partial \overline u}$$
where $V$ is a real-valued homogeneous polynomial of degree $3$ in
$u,\overline u$; with the above expression of $q(u)$, this
corresponds to $c=2\overline{a}$.
 In this case, the following energy is conserved,
$$E=\int _{M}|\nabla u|^2+V(u)\, dx\ .$$
A typical example is
$$V(u)=\frac{1}{2}|u|^2(u+\overline u)\ ,\
q(u)=|u|^2+\frac{1}{2}u^2 \ .$$ This Hamiltonian structure does
not prevent solutions from blowing up in general.
 In the above example, for instance, a purely imaginary constant as Cauchy data leads to a blow up
 solution. However it is possible to give a
classification of all the Hamiltonian quadratic nonlinearities for
which the Cauchy problem associated to (\ref{eq3'}) has a unique
global solution for  small initial data in $H^1_{{\rm
zonal}}({\mathbb{S}}^4)$.

\begin{cor}\label{C1}
Assume $(M,g)$ is the four-dimensional standard sphere and
$c=2\overline{a}$. Then the following assertions are equivalent.\\
i) There exists a subspace $X$ of $C(\mathbb{R}, H^1_{{\rm
zonal}}({\mathbb{S}}^4))$ such that, for every small initial data
 $\|u_0\|_{H^1_{{\rm zonal}}({\mathbb{S}}^4)} \leq \varepsilon$,
 the Cauchy problem (\ref{eq3'}) has a unique global
solution $u \in X$. \\
ii) The parameters $a,b$ satisfy
\begin{equation}\label{Cond}
 \frac{\overline{a}^{2}}{a} = {b}.
\end{equation}
\end{cor}

It would be interesting to know whether blowing up solutions exist
for non small data under property (\ref{Cond}).

 When  property
(\ref{Cond}) is not satisfied, our blowing up solutions are
particularly simple, since they are solutions of the  ordinary
differential equation deduced from (\ref{eq3'}) for
space-independent solutions. Another open problem is of course to
find a wider variety of blowing up solutions for equation
(\ref{eq3'}) in this case.

 \vskip 0.5cm \noindent {\bf Acknowledgements.} This
paper was written while the second author visited the Laboratoire
de Math\'ematiques d'Orsay, supported by the HYKE network. She is
grateful to these two institutions for their hospitality and
support.

\section{Wellposedness via  multilinear estimates}
%Resolution of the nonlinear Schr\"odinger equation using multilinear estimates}

The main step of this section is to prove a result of local
existence in time for initial data in $H^1(M)$ using some
multilinear estimates associated to the nonlinear Schr\"odinger
equation, that we will establish in  Section 3 with a special
attention to the case of the sphere.  For that purpose we follow
closely the ideas of Burq, G\'erard and Tzvetkov (\cite{BGT3},
\cite{BGT8}). In those papers, the authors extended
 to general compact manifolds the  nonlinear methods introduced
by Bourgain (\cite{B1}, \cite{B2}, \cite{B4}) in the context of tori
${\mathbb{R}^d}/{\mathbb{Z}^d}$. Finally,  we achieve the global
wellposedness thanks to the conservation laws.

\subsection{Well-posedness  in Sobolev spaces for the Hartree nonlinearity}

In this subsection we prove that the uniform wellposedness of
(\ref{eq1'}) on $M$ can be deduced from quadrilinear
estimates on solutions of the linear equation. Firstly, we recall
the notion of wellposedness we are going to address.

\begin{definizione}\label{uwp}
Let $s \in \mathbb{R}$. We shall say that the nonlinear Schr\"odinger equation (\ref{eq1'}) is (locally in time) uniformly
well-posed on $H^s(M)$ if, for any bounded subset $B$ of $H^s(M)$, there exists $T>0$ and a Banach space $X_T$
continuously contained into $C([-T,T], H^s(M))$, such that
\begin{description}
  \item[i] For every Cauchy data $u_0\in B$, (\ref{eq1'}) has a unique solution $u \in X_T$.
  \item[ii] If $u_0 \in H^{\sigma}(M)$ for $\sigma > s$, then $u \in C([-T,T], H^{\sigma}(M))$.
  \item[iii] The map $u_0 \in B \mapsto u \in X_T$ is uniformly continuous.
\end{description}
\end{definizione}

The following theorem stresses the general relationship between
uniform wellposedness for equation (\ref{eq1'}) and a certain type
of quadrilinear estimates.

\begin{teo}\label{t3}
Suppose that there exists $C>0$ and $s_0\geq0$  such that for any $f_1, f_2, f_3, f_4 \in
L^2(M)$  satisfying
\begin{equation}\label{Q1}
\mathbf{1}_{\sqrt{1-\Delta} \in [N_j,2N_j]}(f_j) = f_j, \;\;
j=1,2,3,4,
\end{equation}
 one has the following quadrilinear estimates
\begin{equation}\label{QS}
\begin{aligned}
&\sup _{\tau \in {{\mathbb R}}}
\left| \int_{{\mathbb R}} \int_{M}
\chi (t)\,  e^{it \tau}
(1-\Delta)^{-{\alpha}}
 ( u_1 \overline{u}_2) {u_3} \overline{u}_4 dx dt  \right|\\
&\leq C (m(N_1,\cdots ,N_4))^{s_0} \|f_1\|_{L^2(M)} \|f_2\|_{L^2(M)}  \|f_3\|_{L^2(M)}  \|f_4\|_{L^2(M)},\\
&u_j(t)=S(t)f_j,\;\; j=1,2,3,4,\end{aligned}
\end{equation}
where  $\chi \in {\mathcal C}^\infty _0({\mathbb R}) $ is
arbitrary, and $m(N_1,\cdots ,N_4)$ denotes the product of the
smallest two numbers among $N_1, N_2, N_3, N_4$. Then the Cauchy
problem (\ref{eq1'}) is uniformly well-posed in $H^s(M)$ for any
$s > s_0$.
\end{teo}
\begin{proof}
The proof follows essentially the same lines as the one of Theorem
3 in \cite{BGT8} and relies on the use of a suitable class
$X^{s,b} $ of Bourgain-type spaces. We shall sketch it for the
commodity of the reader. We first show that (\ref{QS}) is
equivalent to a quadrilinear estimate in the spaces $X^{s,b} $. We
then prove the crucial nonlinear estimate, from which uniform
wellposedness can be obtained by a contraction argument in
$X_T^{s,b}$. Since this space is continuously embedded in
$C([-T,T], H^s(M))$ provided $b> \frac{1}{2}$, this concludes the
proof of the local well posedness result.

Following the definition in Bourgain \cite{B1} and Burq, G\'erard
and Tzvetkov \cite{BGT3}, we introduce the family of Hilbert
spaces
\begin{equation}
X^{s,b} (\mathbb{R} \times M) = \{v \in \mathcal{S'}(\mathbb{R} \times M) :
(1+|i \partial_t + \Delta|^2 )^{\frac{b}{2}} (1-\Delta)^{\frac{s}{2}}v \in L^2(\mathbb{R} \times M)\}
\end{equation}
for $s, b \in \mathbb{R}$. More precisely, with the notation
$$\langle x\rangle=\sqrt{1+|x|^2}\ ,$$
we have the following definition :

\begin{definizione}
Let $(M,g)$ be a compact Riemannian manifold, and consider the
Laplace operator $-\Delta$ on $M$.
%which is a non negative
%self-adjoint operator with compact resolvent on $L^2(M)$.
 Denote
by $(e_k)$ an $L^2$ orthonormal basis of eigenfunctions of
$-\Delta$, with eigenvalues $\mu_k$, by $\Pi _k$ the orthogonal
projector along $e_k$, and for $s \geq 0$ by $H^s(M)$ the natural
Sobolev space generated by $(I-\Delta)^{\frac{1}{2}}$, equipped
with the following norm
\begin{equation}
\|u\|_{H^s(M)}^2= \sum_k \langle \mu_k \rangle^s \|\Pi _k u
\|_{L^2(M)}^2.
\end{equation}
Then, the space $X^{s,b} (\mathbb{R} \times M)$ is defined as the completion of $C_0^{\infty}(\mathbb{R}_t; H^s(M))$
for the norm
\begin{equation}
\begin{aligned}
&\|u\|_{X^{s,b} (\mathbb{R} \times M)}^2 = \sum_k \|\langle \tau + \mu_k \rangle^b \langle \mu_k \rangle^{\frac{s}{2}}
\widehat{\Pi _k u}(\tau)\|_{L^2(\mathbb{R}_{\tau}; L^2(M))}^2\\
& = \|S(-t)\, u(t,\cdot)\|_{H^b(\mathbb{R}_t; H^s(M))}^2,\\
\end{aligned}
\end{equation}
where $\widehat{\Pi _k u}(\tau)$ denotes the Fourier transform of
$\Pi _k u$ with respect to the time variable.
\end{definizione}
Denoting by $X_T^{s,b}$ the space of restrictions of elements of $X^{s,b}(\mathbb{R} \times M)$ to $]-T,T[ \times M$,
it is easy to prove the embedding
\begin{equation}
\forall b > \frac{1}{2}, \quad X_T^{s,b} \subset C([-T,T], H^s(M)).
\end{equation}
Moreover, we have the elementary property
\begin{equation}
\forall f \in  H^s(M), \quad \forall b >0,
\quad (t,x) \mapsto S(t)f(x) \in  X_T^{s,b}.
\end{equation}

We next reformulate the quadrilinear estimates (\ref{QS}) in the context of $X^{s,b}$ spaces.

\begin{lemma}\label{L1}
Let $s \in \mathbb{R}$. The  following two statements are equivalent:\\
i) For any $f_j \in L^2(M), \; j=1,2,3,4,$ satisfying \eqref{Q1}, estimate (\ref{QS}) holds;\\
ii) For any $b> \frac{1}{2}$ and any $u_j \in
X^{0,b}(\mathbb{R}\times M), \; j=1,2,3,4,$ satisfying
\begin{equation*}
\mathbf{1}_{\sqrt{1-\Delta} \in [N_j,2N_j]}(u_j)=u_j,
\end{equation*}
one has
\begin{equation}\label{15'}
\left| \int_{{\mathbb R}} \int_{M} (1-\Delta)^{-{\alpha}} (u_1
\overline{u}_2) u_3 \overline{u}_4  dx dt  \right| \leq C
(m(N_1,\cdots ,N_4))^{s_0} \prod_{j=1}^4
\|u_j\|_{X^{0,b}(\mathbb{R} \times M)}.
\end{equation}

\begin{proof}
We sketch only the essential steps of the proof of ii) assuming
i), since we follow closely the argument of Lemma 2.3 in
\cite{BGT3}. The reverse implication is easier and will not be
used in this paper.

Suppose first that $u_j$ are supported in time in the interval $(0,1)$
and we select $\chi \in {\mathcal C}^\infty _0({\mathbb
R})$ such that $\chi =1$ on $[0,1]$;
then writing $u^{\sharp}_j(t)=S(-t)u_j(t) $
we have easily
\begin{equation*}
\begin{aligned}
&\left( (1-\Delta)^{-{\alpha}}(u_1 \overline{u}_2) u_3
\overline{u}_4 \right) (t) = {\frac{1}{(2 \pi)^4}} \int_\mathbb{R}
\int_\mathbb{R} \int_\mathbb{R} \int_\mathbb{R}
 e^{it (\tau_1 - \tau_2 +\tau_3 - \tau_4)} \\
 &\times (1-\Delta)^{-{\alpha}}(S(t) \widehat{u}^{\sharp}_1(\tau_1)
  \overline{S(t) \widehat{u}^{\sharp}_2(\tau_2)}) S(t)\widehat{u}^{\sharp}_3(\tau_3)
   \overline{S(t) \widehat{u}^{\sharp}_4(\tau_4)}
  \,d \tau_1 \,d \tau_2  \,d \tau_3 \,d \tau_4,\\
\end{aligned}
\end{equation*}
where $\widehat{u}^{\sharp}_j$ denotes the Fourier transform of ${u}^{\sharp}_j$ with
respect to time. Using i) and the Cauchy-Schwarz inequality in
$(\tau_1, \tau_2, \tau_3, \tau_4)$ (here the assumption $b> \frac{1}{2}$ is used, in
order to get the necessary integrability) yields
\begin{equation*}
\begin{aligned}
\left| \int_{{\mathbb R}\times M} (1-\Delta)^{-{\alpha}} (u_1
\overline{u}_2) u_3 \overline{u}_4  dx dt  \right| &\lesssim
m(N_1,\cdots ,N_4)^{s_0}\, \prod_{j=1}^4 \|
  {\langle \tau\rangle}^b \widehat{u}^{\sharp}_j\|_{L^2(\mathbb{R} \times M)}\\
  & \lesssim m(N_1,\cdots ,N_4)^{s_0} \prod_{j=1}^4 \|u_j\|_{X^{0,b}(\mathbb{R} \times M)}\ .\\
\end{aligned}
\end{equation*}
Finally, by decomposing $u_j(t)= \sum_{n \in \mathbb{Z}} \psi(t-\frac{n}{2}) u_j(t)$
 with a suitable $\psi \in C_0^{\infty}(\mathbb{R})$
supported in $(0,1)$, the general case for $u_j$ follows from the
special case of $u_j$  supported  in the time interval $(0,1)$.
\end{proof}
\end{lemma}

Returning to the proof of Theorem \ref{t3}, there is another way
of estimating the $L^1$ norm of the product
$\left((1-\Delta)^{-{\alpha}} (u_1 \overline{u}_2) u_3
\overline{u}_4 \right)$.
 %Set $P_{N_i}={\bf 1}_{N_i \leq \sqrt {1-\Delta}\leq 2N_i}$\ .

\begin{lemma}\label{interpo}
Assume $\alpha$ as in Theorem 1
and that $u_1, u_2,u_3,u_4$ satisfy
\begin{equation}\label{specloc}
\mathbf{1}_{\sqrt{1-\Delta} \in [N,2N]}  (u_j)=
u_{j}.
\end{equation}
Then, for every $s' > s_0$ there exists $\displaystyle {b'\in
]0,\frac{1}{2}[}$
 such that
\begin{equation}\label{21}
%\begin{aligned}
\left| \int_{{\mathbb R}} \int_{M} (1-\Delta)^{-{\alpha}}
 ( u_1 \overline{u}_2) {u_3} \overline{u}_4 dx dt  \right|
\leq C m(N_1,\cdots ,N_4)^{s'} \prod_{j=1}^4 \|u_j\|_{X^{0,b'}} \
.
%\end{aligned}
\end{equation}
\end{lemma}
\begin{proof}
We split the proof in several steps.

First of all we prove that, for $\alpha> 0$,
\begin{equation}\label{eq15}
\left| \int_{{\mathbb R}} \int_{M}
(1-\Delta)^{-{\alpha}}
(u_1 \overline{u}_2) u_3 \overline{u}_4  dx dt  \right|
 \leq C
m(N_1,\cdots ,N_4)^{2} \prod_{j=1}^4 \| u_j\|_{X^{0,1/4}}.
\end{equation}
By symmetry we have to consider the following three cases:
$$m(N_1,\cdots ,N_4)= N_1 N_2\ , m(N_1,\cdots ,N_4)= N_3 N_4\ ,
\  m(N_1,\cdots ,N_4)= N_1 N_3.$$ In the first case, by a
repeated use of H\"older's inequality, we obtain
\begin{equation*}
\begin{aligned}
&\left| \int_{{\mathbb R}} \int_{M}
(1-\Delta)^{-{\alpha}}
(u_1 \overline{u}_2) u_3 \overline{u}_4  dx dt  \right|\\
& \leq C \|(1-\Delta)^{-{\alpha}}
(u_1 \overline{u}_2)\|_{L^2(\mathbb{R}, L^{\infty }(M))} \| u_3 \overline{u}_4 \|_{L^2(\mathbb{R}, L^{1}(M))},\\
& \leq C \|u_1 \overline{u}_2\|_{L^2(\mathbb{R}, L^{\infty }(M))}
\| u_3 \overline{u}_4 \|_{L^2(\mathbb{R}, L^{1}(M))}\\
& \leq C \|u_1 \|_{L^4(\mathbb{R}, L^{\infty }(M))}
\|{u}_2\|_{L^4(\mathbb{R}, L^{\infty }(M))}
 \| u_3\|_{L^4(\mathbb{R}, L^{2}(M))} \|{u}_4 \|_{L^4(\mathbb{R}, L^{2}(M))},
\end{aligned}
\end{equation*}
where we also used that $(1-\Delta )^{-\alpha}$ is a
pseudodifferential operator of negative order, hence acts on
$L^\infty (M)$. By Sobolev inequality, we infer
\begin{equation*}
\left| \int_{{\mathbb R}} \int_{M} (1-\Delta)^{-{\alpha}} (u_1
\overline{u}_2) u_3 \overline{u}_4  dx dt  \right| \leq C {(N_1
N_2)}^{2} \prod_{j=1}^4\|u_j \|_{L^4(\mathbb{R}, L^{{2}}(M))}\ .
\end{equation*}
By the Sobolev embedding in the time variable for the function
$v(t)=S(-t)u(t)$, we have $X^{0,1/4}\subset {L^4(\mathbb{R},
L^{2}(M))} $,
 and this conclude the proof of the
first case.

In the second case $m(N_1,\cdots ,N_4)= N_3 N_4$ we can proceed in
the same way by writing the integral in the form
\begin{equation*}
\left| \int_{{\mathbb R}} \int_{M}
u_1 \overline{u}_2 (1-\Delta)^{-{\alpha}}(u_3 \overline{u}_4)  dx dt \right|.
\end{equation*}

Finally, when $m(N_1,\cdots ,N_4)= N_1 N_3$, we write the integral
as follows
\begin{equation*}
 \left| \int_{{\mathbb R}} \int_{M}
(1-\Delta)^{-\frac{\alpha}{2}}(u_1 \overline{u}_2) (1-\Delta)^{-\frac{\alpha}{2}}(u_3 \overline{u}_4)  dx dt \right|,
\end{equation*}
and by Cauchy-Schwarz and H\"older's inequalities we estimate it
by
\begin{equation*}\begin{aligned}
&\leq \|(1-\Delta)^{-\frac{\alpha}{2}}(u_1
\overline{u}_2)\|_{L^2(\mathbb{R}, L^{2}(M))}
 \| (1-\Delta)^{-\frac{\alpha}{2}}(u_3 \overline{u}_4)\|_{L^2(\mathbb{R},
 L^{2}(M))}\\
&\leq C \|u_1 \overline{u}_2\|_{L^2(\mathbb{R}, L^{2}(M))}
 \| u_3 \overline{u}_4\|_{L^2(\mathbb{R}, L^{2}(M))}\\
&\leq C  \|u_1 \|_{L^4(\mathbb{R}, L^{{\infty}}(M))}
\|{u}_2\|_{L^4(\mathbb{R}, L^{2}(M))}
 \| u_3\|_{L^4(\mathbb{R}, L^{\infty }(M))} \|{u}_4 \|_{L^4(\mathbb{R},
 L^{2}(M))}\ .
 \end{aligned}
\end{equation*}
Finally we conclude the proof of (\ref{eq15}) by means of
Sobolev's inequality in both space and time variables as above.
\vskip 0.25cm \noindent The second step consists in interpolating
between (\ref{15'}) and (\ref{eq15}) in order to get the estimate
(\ref{21}). To this end we decompose each $u_j$ as follows
\begin{equation*}
u_j= \sum_{K_j} u_{j,K_j}, \qquad  u_{j,K_j}=  \mathbf{1}_{K_j
\leq \langle i \partial_t + \Delta \rangle< 2K_j} (u_j),
\end{equation*}
where $K_j$ denotes the sequence of dyadic integers. Notice that
\begin{equation*}
\|u_j\|^2_{X^{0,b}} \simeq \sum_{K_j}  K_j^{2b} \| u_{j, K_j} \|^2_{L^2(\mathbb{R}\times M)} \simeq
\sum_{K_j}  \| u_{j, K_j} \|^2_{X^{0,b}}.
\end{equation*}
We then write the integral in the left hand side of (\ref{21}) as
a sum of the following elementary integrals,
\begin{equation*}
I(K_1,\cdots ,K_4)=\int_{{\mathbb R}} \int_{M}
(1-\Delta)^{-{\alpha}} (u_{1,K_1} \overline{u}_{2,K_2}) u_{3,K_3}
\overline{u}_{4,K_4} dx dt\ .
\end{equation*}
Using successively (\ref{15'}) and (\ref{eq15}), we estimate these
integrals as
\begin{equation}\label{IK}
|I(K_1,\cdots ,K_4)|\leq C  m(N_1,\cdots ,N_4)^{\sigma } \sum_{
{K_1,K_2,K_3}} (K_1 K_2 K_3K_4)^{\beta }
 \prod_{j=1}^4\|  u_{j, K_j} \|_{L^2}
,\end{equation} where either $(\sigma ,\beta )=(s_0,b)$ for every
$b>1/2$, or $(\sigma ,\beta )=(2,1/4)$. Therefore, for every
$s'>s_0$, there exists $b_1<1/2$ such that (\ref{IK}) holds for
$(\sigma ,\beta )=(s',b_1)$. Choosing $b'\in ]b_1,1/2[$, this
yields
\begin{equation*}\begin{aligned}
&\left| \int_{{\mathbb R}} \int_{M} (1-\Delta)^{-{\alpha}} (u_1
\overline{u}_2) u_3 \overline{u}_4  dx dt  \right|\\ &\leq C\,
m(N_1,\cdots ,N_4)^{s'}\sum _{K_1,\cdots
,K_4}(K_1K_2K_3K_4)^{b_1-b'}\, \prod _{j=1}^4 \| u_j\|
_{X^{0,b'}}\ , \end{aligned}\end{equation*} which completes the
proof, since the right hand side is a convergent series.
\end{proof}

We are finally in position to prove Theorem \ref{t3}.
We can write the solution of the Cauchy problem (\ref{eq1'}) using the Duhamel formula
\begin{equation}\label{D'}
u(t) = S(t) u_0 -i \int_0^t S(t-\tau) \left( (1- \Delta)^{-\alpha}
(|u(\tau)|^2)u(\tau) \right)  d\tau \ .
\end{equation}
 The next lemma
contains the basic linear estimate.
\begin{lemma}\label{3.4'}
Let $b,b'$ such that $0\leq b' < \frac{1}{2}$, $0 \leq b < 1-b'$.
There exists $C>0$ such that, if $T \in [0,1]$, $w(t)= \int_0^t S(t-\tau) f(\tau) d\tau,$
then
\begin{equation}
\|w\|_{ X_T^{s,b}} \leq C T^{1-b-b'} \|f\|_{X_T^{s,-b'}}.
\end{equation}
\end{lemma}
We refer to \cite{G} for a simple proof of this lemma.

The last integral equation (\ref{D'}) can be handled by means of these spaces $X_T^{s,b}$ using  Lemma \ref{3.4'}
as follows
\begin{equation}
\begin{aligned}
\left\| \int_0^t S(t-\tau) \right. & \left.
\left( (1- \Delta)^{-\alpha}  (|u(\tau)|^2)u(\tau) \right)  d\tau \right\|_{X_T^{s,b}} \\
 &\leq C T^{1-b-b'}
\|  \left( (1- \Delta)^{-\alpha}  (|u(\tau)|^2)u(\tau) \right) \|_{X_T^{s,-b'}}.\\
\end{aligned}
\end{equation}
Thus to construct the contraction $\Phi: X^{s,b}_T \rightarrow
X^{s,b}_T,\;\; \Phi(v_i)=u_i, i=1,2$ and to prove the propagation
of regularity ii) in Definition \ref{uwp}, it is enough to prove
the following result.
\begin{lemma}\label{F}
 Let $s > s_0$.
There exists $(b,b') \in \mathbb{R}^2$ satisfying
\begin{equation}\label{42}
0 < b' < \frac{1}{2} <b, \quad b+b' <1,
\end{equation}
and $C>0$ such that for every triple $(u_j), \; j=1,2,3$ in $X^{s,b}(\mathbb{R} \times M)$,
\begin{equation}\label{crucial}
\|   (1- \Delta)^{-\alpha}  (u_1 \overline{u}_2)u_3 \|_{X^{s,-b'}}
\leq C \|u_1\|_{X^{s,b}} \|u_2\|_{X^{s,b}} \|u_3\|_{X^{s,b}}.
\end{equation}
Moreover, for every $\sigma >s$, there exists $C_\sigma $ such
that
\begin{equation}\label{sigma}
\| (1-\Delta )^{-\alpha }(|u|^2)u\| _{X^{\sigma ,-b'}} \leq
C_{\sigma }\| u\| _{X^{s,b}}^2\| u\| _{X^{\sigma ,b}}\ .
\end{equation}
\end{lemma}
\begin{proof}
We only sketch the proof of (\ref{crucial}). The proof of
(\ref{sigma}) is similar. Thanks to a duality argument it is
sufficient to show the following
\begin{equation}\label{B9}
\left| \int_\mathbb{R} \int_{M}  (1- \Delta)^{-\alpha}  (u_1
\overline{u}_2)u_3 \overline{u}_4 dx dt\right| \leq C \left (
\prod _{j=1}^3\|u_j\|_{X^{s,b}}\right ) \|u_4\|_{X^{-s,b'}}.
\end{equation}
%which is equivalent to
%\begin{equation}\label{9}
%\begin{aligned}
%&\left| \int_\mathbb{R} \int_{M} \left( (1- \Delta)^{-\frac{\alpha}{2}}  (u_1 \overline{u}_2)
%(1- \Delta)^{-\frac{\alpha}{2}} (u_3 \overline{u}_4)\right)dx dt\right|\\
%&\leq C \|u_1\|_{X_T^{s,b}} \|u_2\|_{X_T^{s,b}} \|u_3\|_{X_T^{s,b}} \|u_4\|_{X_T^{-s,b'}}.\\
%\end{aligned}
%\end{equation}
The next step is to perform  a  dyadic expansion in the integral
of the left hand-side of (\ref{B9}), this time in the space
variable. We decompose $u_1, u_2,u_3,u_4$ as follows:
\begin{equation*}
u_j=\sum_{N_j} u_{j,N_j}, \qquad   u_{j,N_j}= \mathbf{1}_{\sqrt{1-\Delta} \in [N_j,2N_j]}(u_j).
\end{equation*}
In this decomposition we have
\begin{equation*}
\|u_j\|^2_{X^{s,b}} \simeq \sum_{N_j} N_j^{2s} \|u_{j,N_j}\|^2_{X^{0,b}}
\simeq \sum_{N_j}  \|u_{j,N_j}\|^2_{X^{s,b}}.
\end{equation*}
We introduce now  this decomposition in the left hand side of
(\ref{B9}), and we are left with estimating each term
$$J(N_1,\cdots ,N_4)=\int_\mathbb{R} \int_{M}  (1- \Delta)^{-\alpha}
(u_{1,N_1} \overline{u}_{2,N_2})u_{3,N_3} \overline{u}_{4,N_4} dx
dt\ .$$ Consider the terms with $N_1 \leq N_2 \leq N_3 $ (the
other cases are completely similar by symmetry). Choose $s'$ such
that $s>s'>s_0$. By Lemma \ref{interpo} we can find $b'$ such that
$0 < b' < \frac{1}{2} $ and
\begin{equation}\label{47}
|J(N_1,\cdots ,N_4)| \leq C \sum_{N_j} (N_1 N_2)^{s'}
\prod_{j=1}^4 \|u_{j, N_j}\|_{X^{0,b'}} .
\end{equation}
This is equivalent to
\begin{equation*}
|J(N_1,\cdots ,N_4)| \leq C \sum_{{N_j}} (N_1 N_2)^{s'-s}
\left(\frac{N_4}{ N_3}\right)^{s} \prod _{j=1}^3\| u_{j,N_j}
\|_{X^{s,b'}}  \| { u_{4,N_4}}\|_{X^{-s,b'}}\ .
\end{equation*}
In this series we separate the terms in which $N_4 \leq C N_3$
from the others. For the first ones the series converges thanks to
a simple argument of summation of geometric series and
Cauchy-Schwarz inequality. To perform the summation of the other
terms, it is sufficient to apply the following lemma, which is a
simple variant of Lemma 2.6 in \cite{BGT8}.
\begin{lemma}
Let $\alpha$ a positive number. There exists $C>0$ such that, if
for any $j=1,2,3$, $C\mu_{k_j} \leq \mu_{k_4}$, then for every
$p>0$ there exists $C_p>0$ such that for every $w_j \in L^2(M)$,
$j=1,2,3,4$,
\begin{equation*}
\int_M (1-\Delta)^{-\alpha} (\Pi _{k_1} w_1 \Pi _{k_2} w_2) \Pi
_{k_3} w_3 \Pi _{k_4} w_4 dx \leq C_p\; \mu_{k_4}^{-p}\;
 \prod_{j=1}^4 \|w_j\|_{L^2}.
\end{equation*}
\end{lemma}
\begin{osservazione}
Notice that if $M=\mathbb{S}^4$ the above lemma is trivial since
in that case, by an elementary observation on the degree of the
corresponding sphe\-rical harmonics, we obtain that if $k_4 >
k_1+k_2+k_3$ then the integral (\ref{47}) is zero.
\end{osservazione}
\noindent Finally, the proof of Lemma \ref{F} is achieved by
choosing $b$ such that $\frac{1}{2}<b<1-b'$ and by merely
observing that
$$\|u_j\| _{X^{s,b'}}\leq \|u_j\| _{X^{s,b}}\ ,\ j=1,2,3.$$
\end{proof}

\subsection{Local wellposedness for the quadratic nonlinearity}

\noindent In this subsection, we study the wellposedness theory of
the quadratic nonlinear Schr\"o\-dinger equation posed on $S^4$
\begin{equation}\label{eq3}
{i} \partial_t u + \Delta u =  q(u) ,\;\;\;\;\;
q(u)= a u^2 + b \overline{u}^2+c|u|^2,\\
\end{equation}
with zonal initial data $u(0,x)=u_0(x)$.

\noindent In fact we shall prove Theorem \ref{t2} on every
four-manifold satisfying the trilinear estimates (\ref{TE}). This
is a result of independent interest that we state below.

\begin{teo}\label{t4}
Let $M$ be a Riemannian manifold, let $G$ be a subgroup of
isometries of $M$.  Assuming that there exists $C
>0$ and $s_0$ such that for any $u_1, u_2, u_3 \in
L^2(\mathbb{S}^4)$ $G$-invariant functions on $M$ satisfying
\begin{equation}\label{spectral}
\mathbf{1}_{\sqrt{1-\Delta} \in [N_j,2N_j]}(f_j) = f_j, \;\;
j=1,2,3,
\end{equation}
 one has the trilinear estimates
\begin{equation}\label{TS}
\sup _{\tau \in {{\mathbb R}}}\left| \int_{{\mathbb R}} \int_{M}
\chi (t)\,   e^{it \tau}  \mathcal{T}(u_1, u_2,{u}_3) dx dt
  \right|
 \leq C (\min(N_1,N_2, N_3))^{s_0} \prod_{j=1}^3 \|f_j\|_{L^2} ,
\end{equation}
where $\mathcal{T}(u_1, u_2,{u}_3)=u_1u_2u_3$ or $\mathcal{T}(u_1,
u_2,{u}_3)=u_1u_2\overline {u}_3$ and  $\chi \in {\mathcal
C}^\infty _0({\mathbb R}) $ is arbitrary. Then, for every $s
> s_0$, the Cauchy problem (\ref{eq3})  is uniformly well-posed on
the subspace of $H^s(M)$ which consists of $G$-invariant
functions.
\end{teo}

\begin{proof}
%Accentuare il fatto che qui ci sono stime trilineari nuove che appaiono in maniera naturale su $S^4$.\\
%Tagliare il piu' possibile le cose che gia' stanno su BGT\\
It is close to the one of Theorem \ref{t3} above, so we shall just
survey it. We denote by $L^2_G(M)$, $H^s_G(M)$, $X^{s,b}_G(\mathbb
{R} \times M)$ the subspaces of $L^2(M)$, $H^s(M)$,
$X^{s,b}(\mathbb {R} \times M)$ which consist of $G$-invariant
functions. For the sake of simplicity, we shall focus on the case
$$q(u)=|u|^2+\frac{1}{2}u^2\ .$$
The general case follows from straightforward modifications. As in
the proof of Theorem \ref{t3}, it is enough, for every $s>s_0$, to
show that there exists $b,b'$ such that
$$0<b'<\frac{1}{2}<b<1-b'$$
with the  following estimates,
\begin{equation*}\begin{aligned}
\|u_1 u_2\|_{X^{s,-b'}} \leq C \|u_1\|_{X^{s,b}}
\|u_2\|_{X^{s,b}}\ &,\ \|u_1 \overline{u}_2\|_{X^{s,-b'}} \leq C
\|u_1\|_{X^{s,b}} \|u_2\|_{X^{s,b}}\ ,\\
\| u^2\| _{X^{\sigma ,-b'}}\leq C_{\sigma
}\|u\|_{X^{s,b}}\|u\|_{X^{\sigma ,b}}\ ,\ \| |u|^2&\| _{X^{\sigma
,-b'}}\leq C_{\sigma }\|u\|_{X^{s,b}}\|u\|_{X^{\sigma ,b}}\ ,\
\sigma >s\, ,
\end{aligned}\end{equation*}
where $u_1,u_2,u$ are $G$ -invariant. As before, we focus on the
first set of estimates. Thanks to a duality argument, these
estimates are equivalent to
\begin{equation}\label{T3}\begin{aligned}
\left| \int_{{\mathbb R}} \int_{M} u_1 u_2 \overline{u}_3 dx dt
\right| &\leq C \|u_1\|_{X^{s,b}} \|u_2\|_{X^{s,b}}
\|u_3\|_{X^{-s,b'}}\ ,\\\left| \int_{{\mathbb R}} \int_{M}
\overline {u}_1 u_2 u_3 dx dt \right| &\leq C \|u_1\|_{X^{s,b}}
\|u_2\|_{X^{s,b}} \|u_3\|_{X^{-s,b'}}\ ,
\end{aligned}\end{equation}

 In this way,
%following the same idea of Theorem \ref{t1}
 writing the solution  of the Cauchy problem (\ref{eq3}) using the Duhamel  formula
\begin{equation}\label{D2}
u(t) = S(t) u_0 -i \int_0^t S(t-\tau)
(|u(\tau)|^2+\frac{1}{2}u^2(\tau))\, d\tau,
\end{equation}
and applying Lemma {\ref{3.4'}}, we obtain a contraction  on
$X_T^{s,b}$ proving a result of local existence of the solution to
(\ref{eq3}) on $H^s(M), \; s > s_0 .$ Thus the proof of this
theorem is reduced to establishing the  trilinear estimates
(\ref{T3}) for suitable $s,b,b'$. We just prove the first
inequality in (\ref{T3}). The proof of the second one is similar.

\noindent First we reformulate  trilinear estimates (\ref{TS}) in
the context of Bourgain spaces.

\begin{lemma}
Let $s_0 \in \mathbb{R}$. The  following two statements are equivalent:\\
- For any $f_1, f_2, f_3 \in L^2_G(M)$ satisfying
(\ref{spectral}), estimate (\ref{TS}) holds.
\\
- For any $b> \frac{1}{2}$ and any $u_1,u_2,u_3 \in
X^{0,b}_G(\mathbb{R}\times M)$ satisfying
\begin{equation}\label{specloc2}
\mathbf{1}_{\sqrt{1-\Delta} \in [N_j,2N_j]}(u_j)=u_j, \;\;
j=1,2,3,
\end{equation}
one has
\begin{equation}\label{2.3-D'}
\left| \int_{{\mathbb R}} \int_{M} (u_1 u_2 \overline{u}_3)dx dt
\right| \leq C (\min(N_1, N_2, N_3))^{s_0} \prod
_{j=1}^3\|u_j\|_{X^{0,b}}
%(\mathbb{R} \times M)}
.
\end{equation}
\end{lemma}
\begin{proof}
The proof of this lemma follows   lines of Lemma \ref{L1} above.
First we assume that $u_1,u_2,u_3$ are supported for $t\in [0,1]$,
and we select $\chi \in {\mathcal C}^\infty _0({\mathbb R})$ such
that $\chi =1$ on $[0,1]$. We set $u^{\sharp}_j(t)=S(-t)u_j(t) $.
Using the Fourier transform, we can write
%\begin{equation}
%f(t) = \frac{1}{2\pi} \int_{-\infty}^{\infty} e^{it \tau_1} e^{it \Delta} \widehat{F}(\tau_1), \;\;
%g(t) = \frac{1}{2\pi} \int_{-\infty}^{\infty} e^{it \tau_2} e^{it \Delta} \widehat{G}(\tau_2), \;\;
%h(t) = \frac{1}{2\pi} \int_{-\infty}^{\infty} e^{it \tau_3} e^{it \Delta} \widehat{H}(\tau_3)
%\end{equation}
%and hence
\begin{equation*}
\begin{aligned}
&\left| \int_{{\mathbb R}}\int_{M}  u_1 u_2\overline{u}_3 dx dt \right|\\
%&{{(2 \pi})^{-3}}
% &\left| \int_{{\mathbb R}} \int_{M} \int_\mathbb{R} \int_\mathbb{R}
%\int_\mathbb{R}
 %e^{it (\tau_1 + \tau_2 -\tau_3)} (S(t) \widehat{u}^{\sharp}_1(\tau_1) S(t) \widehat{u}^{\sharp}_2(\tau_2)
  %\overline{S(t) \widehat{u}^{\sharp}_3(\tau_3)})
  %\,d \tau_1 \,d \tau_2  \,d \tau_3 dx dt \right| \\
& \leq C \int_{\tau_1} \int_{\tau_2} \int_{\tau_3} \left| \int_\mathbb{R}
 \int_{M}\chi (t)  e^{it \tau}
\prod _{j=1}^3S(t)\widehat{u}^{\sharp}_j(\tau_j)
 dx dt \right| d\tau_1 d\tau_2 d\tau_3,\\
\end{aligned}
\end{equation*}
where $\tau=(\tau_1+\tau_2-\tau_3)$. Supposing for instance $N_1
\leq N_2 \leq N_3$ and applying (\ref{TS}) we obtain that the
right hand side is bounded by
\begin{equation*}
\leq C N_{1}^{s_0} \int_{-\infty}^{\infty} \int_{-\infty}^{\infty}
\int_{-\infty}^{\infty}
\|\widehat{u}^{\sharp}_1(\tau_1)\|_{L^2(M)}
\|\widehat{u}^{\sharp}_2(\tau_2)\|_{L^2(M)}
\|{\widehat{u}^{\sharp}_3}(\tau_3)\|_{L^2(M)} d\tau_1 d\tau_2
d\tau_3.
\end{equation*}
We conclude the proof as in the proof of Lemma \ref{L1} in section
2, using the Cauchy-Schwarz inequality in $(\tau_1,\tau_2,
\tau_3)$, and finally decomposing each $u_j$ by means of the
partition of unity
$$1=\sum _{n\in {\mathbb Z}}\psi \left (t-\frac{n}{2}\right )\ ,$$
where $\psi \in {\mathcal C}^\infty _0([0,1])$.
\end{proof}

\begin{lemma}\label{interpo2}
For every $s' > s_0$ there exist $b'$
 such that $0 < b' < \frac{1}{2}$ and, for every $G$-invariant
 functions $u_1, u_2,u_3$ satisfying (\ref{specloc2}),
\begin{equation}
%\begin{aligned}
\left| \int_{{\mathbb R}} \int_{M}
 ( u_1 {u}_2  \overline{u}_3) dx dt  \right|
\leq C \min(N_1,N_2,N_3)^{s'} \prod _{j=1}^3 \|u_j\|_{X^{0,b'}} \
.
%\end{aligned}
\end{equation}
\end{lemma}
\begin{proof}
Following the same lines of the proof of Lemma \ref{interpo}, it
is enough to establish
\begin{equation}\label{9'}
\left| \int_{{\mathbb R}} \int_{M}  ( u_1 {u}_2  \overline{u}_3)
dx dt \right| \leq C \min(N_1, N_2, N_3)^{2} \prod_{j=1}^3
\|u_j\|_{X^{0,\frac{1}{6}}(\mathbb{R} \times M)}.
\end{equation}
Then the lemma follows by interpolation with (\ref{2.3-D'}).
Indeed, assuming for instance $N_1\leq N_2 \leq N_3$, we apply the
H\"older inequality as follows,
\begin{equation*}
\left| \int_{{\mathbb R}} \int_{M} ( u_1 {u}_2  \overline{u}_3) dx dt \right| \leq C
\|u_1\|_{L^3({{\mathbb R}}, L^{\infty}(M))} \|u_2\|_{L^3({{\mathbb R}}, L^{2}(M))}
\|u_3\|_{L^3({{\mathbb R}}, L^{2}(M))}
\end{equation*}
and using the Sobolev embedding we obtain
\begin{equation*}
\leq C (N_1)^{{2}}\|u_1\|_{L^3({{\mathbb R}}, L^{2}(M))}
\|u_2\|_{L^3({{\mathbb R}}, L^{2}(M))} \|u_3\|_{L^3({{\mathbb R}}, L^{2}(M))}.
\end{equation*}
By  the Sobolev embedding in the time variable for function
$v(t)=S(-t)u(t)$, we know that
\begin{equation*}
\|u\|_{L^{3}({\mathbb{R}}, L^2(M))} \leq
\|u\|_{X^{0,\frac{1}{6}}({\mathbb{R}} \times M)}\,
\end{equation*}
and this completes the proof.

\end{proof}

Let us sketch the last part of the proof of Theorem\ref{t4}. We
decompose $u_1, u_2,u_3$ as follows:
\begin{equation*}
u_j=\sum_{N_j} u_{j,N_j}, \qquad   u_{j,N_j}= \mathbf{1}_{\sqrt{1-\Delta} \in [N_j,2N_j]}(u_j).
\end{equation*}

We introduce   this decomposition in the left hand side of
(\ref{T3}) and we use Lemma \ref{interpo2}. Supposing now for
simplicity that $N_1 \leq N_2 $, we obtain that for any $s'
> s_0$ we can find $b'$ such that $0 < b' < \frac{1}{2} $  and
\begin{equation}
\left| \int_\mathbb{R} \int_{M}   u_1 {u}_2 \overline{u}_3 dx
dt\right| \leq C \sum_{{N_j}} (N_1)^{s'-s} \left(\frac{N_3}{
N_2}\right)^{s} \|  u_1 \|_{X^{s,b'}} \|  u_2\|_{X^{s,b'}} \| {
u_3}\|_{X^{-s,b'}}
\end{equation}
for any $s> s' > s_0$. Notice that the summation over $N_1$ can be
performed via a crude argument of summation of geometric
 series.
As for the summation over $ N_2,N_3$, following the same proof as
in Section 2.1, we conclude by observing that the main part of the
series corresponds to the constraint $N_3 \lesssim N_2$.
\end{proof}

\subsection{Conservation laws and global existence for the Hartree nonlinearity}

Next we prove that for an initial datum $u_{0}\in H^{1}(M)$, the
local solution
 of the Cauchy problem (\ref{eq1'})
obtained above can be extended to a global solution $u \in
C(\mathbb{R}, H^1(M))$.

By the definition of uniform wellposedness, the lifespan $T$ of
the local solution $u\in C([0,T), H^1(M))$ depends only on the
$H^{1}$ norm of the initial datum. Thus, in order to prove that
the solution can be extended to a global one, it is sufficient to
show that the $H^{1}$ norm of $u$ remains bounded on any finite
interval $[0,T)$. This is a consequence of the following
conservation laws, which can be proved by means of the multipliers
$\overline u$ and $\overline u_t$,
\begin{equation}\begin{aligned}
&\int _{M}|u(t,x)|^2\, dx= Q_0\ ;\\\ &\int _{M}|\nabla
u(t,x)|_g^2+\frac{1}{2}|(1-\Delta )^{-\alpha /2}(|u|^2)(t,x)|^2\,
dx=E_0\ .
\end{aligned}\end{equation}

\begin{osservazione}
Notice that a similar argument can be applied in the case of an
attractive Hartree nonlinearity,  at least when $\alpha>1$.
Indeed, consider the focusing Schr\"odinger equation
\begin{equation*}
    iu_{t}+\Delta u=-(1-\Delta)^{-\alpha}(|u|^{2})u,
\end{equation*}
where the nonlinear term has the opposite sign.
Computing as above, we obtain the conservation of energy
\begin{equation*}
    \|\nabla  u\|^{2}_{L^{2}(M)}
    -\frac{1}{2}\|(1-\Delta)^{-\alpha/2}(|u|^{2})\|^{2}_{L^{2}}=const,
\end{equation*}
but now the energy $E(t)$ does not control the $H^{1}$ norm
of $u$. However, we can write
\begin{equation*}
    \|\nabla  u\|^{2}_{L^{2}}\leq
    C+C\|(1-\Delta)^{-\alpha/2}(|u|^{2})\|^{2}_{L^{2}},
\end{equation*}
and by Sobolev embedding we have
\begin{equation*}
    \|(1-\Delta)^{-\alpha/2}(|u|^{2})\|^{2}_{L^{2}}\leq
    C\|\;|u|^{2}\|_{L^{q}}^{2}\equiv
    C\|u\|_{L^{2q}}^{4},\qquad
    \frac1q=\frac12+\frac\alpha4,
\end{equation*}
so that we obtain, with $p=2q$,
\begin{equation*}
    \|\nabla  u\|_{L^{2}}\leq
    C+C\|u\|_{L^{p}}^{2},\qquad
    \frac1p=\frac14+\frac\alpha8.
\end{equation*}
We now use the Gagliardo-Nirenberg inequality (for $d=4$)
\begin{equation*}
    \|w\|_{L^{p}}^{p}\leq
    C(\|w\|_{L^{2}}^{p-(p-2)\frac d2}
    \|\nabla  w\|_{L^{2}}^{(p-2)\frac
    d2}+\|w\|_{L^{2}}^{p})
\end{equation*}
and we obtain
\begin{equation*}
    \|\nabla  u\|_{L^{2}}\leq
    C(1+\|u\|_{L^{2}}^{2})+
    C\|u\|_{L^{2}}^{2-4(p-2)/p}
    \|\nabla  u\|_{L^{2}}^{4(p-2)/p}.
\end{equation*}
Notice that, as in the defocusing case above,  the $L^{2}$ norm of
$u$ is a conserved quantity. If the power $4(p-2)/p$ is strictly
smaller than 1, we infer that the $H^{1}$ norm of $u$ must remain
bounded. In other words, we have proved global existence provided
\begin{equation*}
    4\cdot\frac{p-2}p<1\quad\iff\quad \alpha>1.
\end{equation*}
\end{osservazione}
\end{proof}

\subsection{Studying the global existence for the quadratic nonlinearity}

\begin{proposizione}
Let $(M,g)$ be  a four-dimensional Riemannian manifold satisfying
the assumptions of Theorem \ref{t4}.  There exists $\varepsilon
>0$ and a subspace $X$ of $C(\mathbb{R}, H^1_{G}(M))$ such
that, for every initial data $u_0\in H^1_G(M)$ satisfying
 $\|u_0\|_{H^1} \leq \varepsilon$,
 the Cauchy problem (\ref{eq3'}), where $q(u)= (\mathrm{Re}\, u)^2$, has a unique global
solution $u \in X$.
\end{proposizione}
\begin{proof}
By Theorem \ref{t4}, we obtain that for an initial datum $u_{0}\in
H^{1}_{G}(M)$, there exists a local solution
 of the Cauchy problem
 \begin{equation*}
\begin{cases}
&{i} \partial_t u + \Delta u = \,(\mathrm{Re}\, u)^2 ,\\
&u(0,x)=u_0(x).
\end{cases}
\end{equation*}
By the definition of  uniform wellposedness, the lifespan $T$ of
the local solution $u\in C([0,T), H^{1}_G(M))$ only depends on a
bound of the $H^{1}$ norm of the initial datum. Thus, in order to
prove that the solution can be extended to a global one, it is
sufficient to show that the $H^{1}$ norm of $u$ remains bounded on
any finite interval $[0,T)$. This is a consequence of the
following conservation laws and of a suitable assumption of
smallness on the initial data. Notice that
\begin{equation*}
\partial_t \left( \int_M u(t,x)\, dx \right) = -i   \int_{M} (\mathrm{Re}\, u)^2\, dx,
\end{equation*}
from which
\begin{equation}\label{C}
  \int_M {\mathrm{Re}}\, u(t,x)\, dx  = const.
\end{equation}
Moreover the following energy is conserved,
\begin{equation}
 \int _M|\nabla
u(t,x)|^2+\frac23  (\mathrm{Re}\, u(t,x))^3 \, dx=E_0\ .
\end{equation}
 Consequently  we can write
\begin{equation*}
    \|\nabla  u\|^{2}_{L^{2}}\leq
    E_0 +C \left|\int _M  (\mathrm{Re}\, u )^3 \right|.
\end{equation*}
Since by Gagliardo-Nirenberg inequality we have
\begin{equation*}
 \left|\int _M  (\mathrm{Re}\, u )^3 dx\right| \leq
 C \|\mathrm{Re} \, u\|_{L^2}  \| \nabla (\mathrm{Re} \, u) \|_{L^2}^2 +  \|(\mathrm{Re}\, u )\|_{L^2}^3,
\end{equation*}
%\begin{equation*}
 %\left|\int _M  (\mathrm{Re}\, u )^3 dx\right| \leq \left(\int _M |u|^2 \,dx\right)^{\frac12} \|u\|_{H^1}^2,
%\end{equation*}
and by the following inequality
%(Poincare' inequality on $Reu$)
%\begin{equation*}
%\left(\int _M |u|^2 dx\right) \leq C \left|\int _M u \,  dx\right|^2 + \|\nabla u\|_{L^2}^2,
%\end{equation*}
\begin{equation*}
\|\mathrm{Re} \, u\|_{L^2} \leq C \left|\int_M \mathrm{Re} \, u
\,dx \right| + \|\nabla (\mathrm{Re}\, u)\|_{L^2},
\end{equation*}
we deduce that
\begin{equation*}
    \|\nabla  u\|^{2}_{L^{2}}\leq
    E_0 + C \left( \left|\int _M \mathrm{Re}\, u  \,dx\right| + \|\nabla u\|_{L^2}\right) \|\nabla u\|_{L^2}^2.
\end{equation*}
Thanks to (\ref{C}) we know that
\begin{equation*}
\left|\int_M \mathrm{Re} \, u\, dx \right| \leq \|u_0\|_{L^1(M)}
\leq C \|u_0\|_{H^1{(M})},
\end{equation*}
thus we obtain
\begin{equation*}
    \|\nabla  u\|^{2}_{L^{2}}\leq
    E_0 + C \left(  \|u_0\|_{H^1}+ \|\nabla u\|_{L^2}\right) \|\nabla u\|_{L^2}^2.
\end{equation*}
%If initial data $  \|u_0\|_{H^1{({\MATHBB{S}}^4})} \leq \varepsilon$
%\begin{equation*}
 %   \|\nabla  u\|^{2}_{L^{2}}\leq
  %C  E_0 + C   \|\nabla u\|_{L^2}^3.
%\end{equation*}
Assuming that
\begin{equation*}
  \|  u_0\|_{H^1} \leq \varepsilon,
\end{equation*}
we infer, by a classical bootstrap argument, that $\| \nabla u\| $
cannot blow up, as well as $\| {\mathrm Re}\, u\| _{L^2}$. Using
again the evolution law of the integral of $u$, this implies that
this integral cannot blow up, and completes  the proof of the
proposition.
\end{proof}

Notice that the  proof above extends without difficulty to
$q(u)=c(\mathrm{Re}\;u)^2$, for any real number $c$. If $(M,g)$
satisfies the assumptions of Theorem \ref{t4}, we can now prove
that the conclusions of Corollary \ref{C1} hold on $M$.
\begin{proof}

Let $q(u)=a u^2 + b \overline{u}^2+ 2 \overline{a}|u|^2$ . The
idea is to transform the equation into an equivalent one using the
change of unknown $u=\omega v$, with $|\omega|=1$, and then impose
conditions on $a,b$ such that the transformed equation is of the
special type corresponding to $q(u)=c(\mathrm{Re}\;u)^{2}$  for
which, thanks to Proposition 1, we know that the solution is
global. Thus we try to impose
\begin{equation*}
    q(\omega v)=c\omega (\mathrm{Re}\;v)^{2}
\end{equation*}
for some $c\in\mathbb{R}$ and some $\omega $ with $|\omega|=1$,
and we obtain the polynomial identity
\begin{equation*}
    a\omega ^2v^2+
 + b \overline{\omega}^{2}\overline {v}^{2}+
    2\overline {a}|v|^2=
    \frac {c \omega} 4(v+ \overline {v})^{2}.
\end{equation*}
Equating the coefficients of the two polynomials we obtain
\begin{equation*}
a=c\frac{\overline \omega }{4},\qquad b=c\frac{\omega ^3}{4}
\end{equation*}
and this is equivalent to
\begin{equation*}
\frac{\overline{a}^{2}}{a} = { b}\ .
\end{equation*}

\noindent Conversely, we prove that if this condition is not
satisfied, it is always possible to construct small energy
solutions which blow up in a finite time. We take as initial datum
a constant in the form
\begin{equation*}
    u_{0}(x)=\omega y_{0},\qquad
       y_{0}\in\mathbb{R}\setminus \{ 0\}\ ,\quad|\omega|=1,
\end{equation*}
and then the equation
reduces to the ordinary differential equation
\begin{equation*}
    iu_{t}=q(u),\qquad u(0)=\omega y_{0}.
\end{equation*}
Defining $y(t)=u(t)/\omega$, we see that $y(t)$ is a solution
of the equation
\begin{equation*}
    i\omega y'(t)=q(u)=y^{2}q(\omega)
\end{equation*}
which can be written
\begin{equation*}
    y'(t)=-iq(\omega)\overline{\omega}\;y^{2},\qquad
    y(0)=y_{0}\in\mathbb{R}
\end{equation*}
The solution can be written explicitly as
\begin{equation*}
    y(t)=\frac 1{y_{0}^{-1}+iq(\omega)\overline{\omega}t}
\end{equation*}
and is not global if and only if $q(\omega)\overline{\omega}$ is purely imaginary.
Thus to conclude the proof it is sufficient to show that we can
find an $\omega$ such that
\begin{equation*}
    q(\omega)\overline{\omega}\equiv
      a\omega+b\overline{\omega}^{3}+2\overline a\overline\omega
      \quad\text{is purely imaginary (and not 0)}.
\end{equation*}
Writing $a=Ae^{i\alpha}$, $b=Be^{i\beta}$, $\omega=e^{i\theta}$
with $A,B\geq0$, this is equivalent to finding a simple zero for
the following function
\begin{equation*}
    f(\theta) = 3A\cos(\alpha+\theta)+B\cos(\beta-3\theta).
\end{equation*}
Observe that the average of $f$ vanishes. A point where the sign
of $f$ changes cannot be a double zero unless it is a triple zero,
and a straightforward calculation shows that this corresponds
exactly to the case $A=B$ and $3\alpha+ \beta= 2k\pi$, namely
$\frac{\overline{a}^{2}}{a} = {b}$. Hence, if this condition is
not satisfied, $f$ has a simple zero. This completes the proof.

\end{proof}

\section{Multilinear estimates}
In this section  we establish  multilinear estimates, which,
combined with Theorems \ref{t3} and \ref{t4}, yield Theorems
\ref{t1'} and \ref{t2}. We recall that $S(t)=e^{it\Delta }\ .$

\subsection{Quadrilinear estimates}
 This subsection is devoted to the proof of quadrilinear estimates (\ref{QS})
 with $s_0<1$ on arbitrary four-manifolds with $\alpha >1/2$, and
 on the sphere $\mathbb {S}^4$ with $\alpha >0$. In view of
 subsections 2.1 and 2.3, this will complete the proof of Theorem \ref{t1'}.

\begin{lemma}\label{BSS}
Let $\alpha > \frac{1}{2}$, $s_0 = \left( \frac{3}{2} - \alpha
\right)$  and let $(M,g)$ a compact four-dimensional Riemannian
manifold. Then there exists $C>0$ such that for any $f_1, f_2 \in
L^2(M)$ satisfying
\begin{equation}\label{}
\mathbf{1}_{\sqrt{1-\Delta} \in [N,2N]}(f_1)=f_1, \;\;\;
\mathbf{1}_{\sqrt{1-\Delta} \in [L,2L]}(f_2)=f_2,
\end{equation}
one has the following bilinear estimate:
\begin{equation}\label{bill}
\|(1-\Delta)^{-\frac{\alpha}{2}} (u_1  u_2) \|_{L^2((0,1)\times
M)} \leq C (\min(N,L))^{s_0} \|f_1\|_{L^2(M)} \|f_2\|_{L^2(M)}\ ,
\end{equation}
with $u_j(t)=S(t)f_j$.
\end{lemma}
\begin{proof}
By symmetry, it is not restrictive to assume that $N \leq L$. The
Sobolev embedding implies
\begin{equation*}
\|(1-\Delta)^{-\frac{\alpha}{2}} (u_1  u_2)\|_{L^2((0,1) \times M)}
\leq C \|u_1  u_2\|_{L^2((0,1), L^{q}(M))}, \;\; \frac{1}{q}=
\frac{1}{2} + \frac{\alpha}{4},
\end{equation*}
and applying the Hold\"er inequality we obtain
\begin{equation*}
\|(1-\Delta)^{-\frac{\alpha}{2}} (u_1  u_2)\|_{L^2((0,1) \times
M)}\leq C \| u_1 \|_{L^{2}((0,1), L^{\frac{4}{\alpha}}(M))} \|
u_2\|_{L^{\infty}((0,1), L^{2}(M))}.
\end{equation*}
Thanks to  the conservation of the $L^{2}$ norm we can bound the
last factor with the $L^{2}$ norm of $f_2$; on the other hand, the
$L^{2}L^{4/\alpha}$ term can be bounded using the Strichartz
inequality on compact manifolds established by Burq, G\'erard,
Tzvetkov in \cite{BGT4}  (see Theorem 1), which reads, in this
particular case,
$$\| u_1 \|_{L^{2}((0,1), L^4(M))}\leq C\,
N^{1/2}\| f_1\| _{L^2(M)}\ .$$
 Combining this estimate  with the Sobolev
inequality, we obtain \eqref{bill} as claimed.
\end{proof}

\begin{proposizione}\label{P1}
Let $\alpha > \frac{1}{2}$, $s_0 > \left( \frac{3}{2} - \alpha
\right)$  and let $(M,g)$ a compact four dimensional Riemannian
manifold. Then there exists $C>0$ such that for any $f_1, f_2,
f_3, f_4 \in L^2(M)$  satisfying
\begin{equation*}
\mathbf{1}_{\sqrt{1-\Delta} \in [N_j,2N_j]}(f_j) = f_j, \;\;
j=1,2,3,4,
\end{equation*}
one has the following quadrilinear estimate for $u_j(t)=S(t)f_j$:
\begin{equation}\label{QSS}
\begin{aligned}
&\sup _{\tau \in {{\mathbb R}}}
\left| \int_{{\mathbb R}} \int_{M}
\chi (t)\,  e^{it \tau}
(1-\Delta)^{-{\alpha}}
 ( u_1 \overline{u}_2) {u_3} \overline{u}_4 dx dt  \right|\\
&\leq C (m(N_1,\cdots ,N_4))^{s_0} \|f_1\|_{L^2(M)} \|f_2\|_{L^2(M)}  \|f_3\|_{L^2(M)}  \|f_4\|_{L^2(M)},\\
\end{aligned}
\end{equation}
where  $\chi \in {\mathcal C}^\infty _0({\mathbb R}) $ is
arbitrary and $m(N_1,\cdots ,N_4)$ is the product of the smallest
two numbers among $N_1, N_2, N_3, N_4$.
\end{proposizione}
\begin{proof}
The proof of our quadrilinear estimate (\ref{QSS}) when
$m(N_1,\cdots ,N_4)= N_1 N_3$ follows directly by the
Cauchy-Schwarz inequality and Lemma \ref{BSS}. In fact, assuming
for instance that $\chi$ is supported into $[0,1]$, we have
\begin{equation*}
\begin{aligned}
&I \equiv \sup _{\tau \in {{\mathbb R}}}
\left| \int_{{\mathbb R}} \int_{M}
\chi (t)\,  e^{it \tau}
(1-\Delta)^{-{\alpha}}
(u_1 \overline{u}_2){u}_3 \overline{u}_4 dx dt  \right|\\
& \leq C
\|(1-\Delta)^{-\frac{\alpha}{2}} (u_1 \overline{u}_2)\|_{L^2((0,1) \times M)}
\|(1-\Delta)^{-\frac{\alpha}{2}} ({u}_3 \overline{u}_4)\|_{L^2((0,1) \times M)}\\
&\leq C (m(N_1,\cdots ,N_4))^{s_0} \|f_1\|_{L^2(M)}
\|f_2\|_{L^2(M)} \|f_3\|_{L^2(M)}  \|f_4\|_{L^2(M)},\end{aligned}
\end{equation*}
by applying (\ref{bill}). By symmetry, it remains to consider only
the case $$m(N_1,\cdots ,N_4)= N_1 N_2\ .$$ By the
self-adjointness of $(1-\Delta)$, H\"older's inequality and
Sobolev's inequality we have
\begin{equation*}\begin{aligned}
I &\leq C \|u_1 \overline{u}_2\|_{L^1((0,1), L^{q'}(M))}
\|(1-\Delta)^{-\alpha} (u_3 \overline{u}_4)\|_{L^\infty((0,1),
L^{q}(M))}\\
 &\leq C \|u_1 \overline{u}_2\|_{L^1((0,1), L^{q'}(M))}
\|u_3 \overline{u}_4\|_{L^\infty((0,1), L^{1}(M))},
\end{aligned}
\end{equation*}
provided $\frac{1}{q} > 1- \frac{\alpha}{2}$. Using again
H\"older's inequality, we infer
\begin{equation*}
I \leq C \prod _{j=1,2}\|u_j\|_{L^2((0,1), L^{2q'}(M))}
\prod_{k=3,4} \|u_k\|_{L^\infty((0,1), L^{2}(M))} \ .
\end{equation*}
Conservation of energy implies that $ \| {u}_k\|_{L^\infty((0,1),
L^{2}(M))}= \|f_k\|_{L^2(M)}$. On the other hand by Sobolev
embedding we have
\begin{equation*}
 \|u_j\|_{L^2((0,1), L^{2q'}(M))} \leq C N_j^{\frac{2}{q}-1}
 \|u_j\|_{L^2((0,1), L^{4}(M))} .
\end{equation*}
Now we can apply the above-mentioned Strichartz estimate of
\cite{BGT4} to obtain
\begin{equation*}
 \|u_j\|_{L^2((0,1), L^{2q'}(M))}\leq C N_j^{\frac{2}{q}-\frac12}
 \|f_j\|_{ L^{2}(M))}\ .
\end{equation*}
Since $$s_0={\frac{2}{q}-\frac12} > \frac{3}{2}-\alpha,$$ and
$s_0$ can be arbitrarily close to $\frac{3}{2}-\alpha$, the proof
 is complete.
\end{proof}
\begin{osservazione}
In this case, an iteration scheme for solving can be performed as
in \cite{BGT4}, avoiding the use of Bourgain spaces, making in
$X_T= C([0,T], H^1) \cap L^2([0,T], H^\sigma_4)$ .
\end{osservazione}

\noindent On the four dimensional sphere, endowed with its
standard metric, the precise knowledge of the spectrum  $\mu_k=
k(k+3), \; k \in \mathbb{N}$ makes it possible to improve
 our quadrilinear estimate. We proceed in several steps, starting
 with an estimate on the product of two spherical harmonics.

\begin{lemma}\label{BB}
Let $\alpha\in ]0,\frac{1}{2}]$ and let $s_0 =1- \frac{3
\alpha}{4}$.
 There exists $C>0$  such that for any $H_{n}, \widetilde{H}_l$ spherical harmonics on $\mathbb{S}^4$
of degree $n,l$ respectively,
the following bilinear estimate holds:
\begin{equation}\label{eqBB}
\|(1-\Delta)^{-\frac{\alpha}{2}} (H_{n} \widetilde{H}_{l})
\|_{L^2( \mathbb{S}^4)} \leq C (1+\min({(n,l}))^{s_0}
\|H_{n}\|_{L^2(\mathbb{S}^4)}
\|\widetilde{H}_{l}\|_{L^2(\mathbb{S}^4)}.
\end{equation}
\end{lemma}
\begin{proof}
It is not restrictive to assume that $1\leq n \leq l$. We shall
adapt the proof of multilinear estimates in
\cite{BGT6},\cite{BGT7}, using the approach described in
\cite{BGT3}.

Writing $$h=(n(n+3))^{-1/2}\ ,\ \widetilde h=(l(l+3)^{-1/2}\ ,$$
the equations satisfied by the eigenfunctions $H_{n}, \widetilde
H_l$ read
$$h^{2}\Delta H_n+H_n=0\ ,\  \tilde h^2\Delta \widetilde H_l+\widetilde H_l=0\ .$$
In local coordinates, these are semiclassical equations, with
principal symbol
$$p(x,\xi )=1-g_x(\xi ,\xi )\ .$$
We now decompose $H_{n}$ and $H_{l}$ using a microlocal partition
of unity  with semi-classical cut-off of the form $\chi(x,hD)$,
$\widetilde \chi (x,\widetilde hD)$ respectively.  When
$$\mathrm{supp}\chi(x,\xi) \cap \{g_x(\xi,\xi)=1\}= \emptyset,$$
i.e. in the "elliptic" case, the estimates are quite strong : we
have, for all $s$, $p$,
\begin{equation}\label{chi}
\||D_x|^s \chi(x, hD_x) H_n\|_{L^2(S^4)} \leq C_{s,p}\,  h^p
\|H_n\|_{L^2(S^4)},
\end{equation}
with similar estimates for $\tilde H_l$. Consequently, it is
sufficient to estimate
\begin{equation}
\|(1-\Delta)^{-\frac{\alpha}{2}} ({\chi}(x, hD_x)H_{n} \,
\widetilde{\chi}(x, \widetilde h
D_x)\widetilde{H}_{l})\|_{L^2(S^4)}
\end{equation}
when  cut-off functions $\chi ,\widetilde \chi $ are localized
near the characteristic set $$\{g_x(\xi,\xi)=1\}\ .$$ Refining the
partition of unity, we may assume that the supports of $\chi $,
$\widetilde \chi $ are contained in small neighborhoods of
$(m,\omega )$, $(m,\widetilde \omega )$ where $m\in M$ and $\omega
,\widetilde \omega $ are covectors such that
$$g_m(\omega ,\omega )=g_m(\widetilde \omega ,\widetilde \omega
)=1\ .$$ Notice that functions $u={\chi}(x, hD_x)H_{n}\ ,\
\widetilde u= \widetilde{\chi}(x, \widetilde h
D_x)\widetilde{H}_{l}$ are compactly supported and satisfy
$$p^w(x,hD)u=hF\ ,\ p^w(x,\widetilde hD)\widetilde u=\widetilde
h\widetilde F\ ,$$ where $\| F\| _{L^2}\lesssim \| H_n\| _{L^2}$
and $\| \widetilde F\| _{L^2}\lesssim \| \widetilde H_l\| _{L^2}$.

Set $g_x(x,\xi)=\langle A(x)\xi,\xi\rangle$. Choose any system
$(x_1, \ldots, x_4)$ of linear coordinates on $\mathbb{R}^4$ such
that
$$\langle A(m)\omega,  dx_1 \rangle \neq 0\quad  {\rm and}\quad
\langle A(m)\widetilde \omega, dx_1 \rangle \neq 0\ .$$
%The next step is to remark that
%after this localization the symbol of $p(x,\xi)$ (and of $\tilde p(x,\xi)$)
%can be factored as
Then, on the supports of $\chi $ and $\widetilde \chi$, one can
factorize the symbol of the equation as
$$p(x,\xi) = e(x,\xi)(\xi_{1}-q(x,\xi'))\ ,\ p(x,\xi) = \widetilde e(x,\xi)(\xi_{1}-\widetilde q(x,\xi')),$$
%for some $j=1,\dots,n$ (we can take e.g. $j=1$),
where $e,\widetilde e$ are  elliptic symbol while $q,\widetilde q$
are real valued symbols. In other words, we can reduce the
equations for $u,\widetilde u$ to  \emph{evolution} equations with
respect to the variable $x_{1}$. Notice that
$\xi'\in\mathbb{R}^{d-1}=\mathbb{R}^{3}$, i.e., the spatial
dimension of these evolution equations is $3$. Moreover, since the
second fundamental form of the characteristic ellipsoid $\{ \xi
:g_m(\xi ,\xi )=1\} $ is non degenerate, the Hessian of
$q,\widetilde q$ with respect to the $\xi '$ variables does not
vanish on the supports of $\chi ,\widetilde \chi $ respectively.

Therefore we can apply to this equation the (local)
three-dimensional Strichartz estimates (see Corollary 2.2 of
\cite{BGT3}  for more details). We conclude
 that $u$ satisfies the $3$-dimensional
 semiclassical Strichartz estimates in the following form:
\begin{equation}\label{eq.strch}
    \|u\|_{L^{p}_{x_{1}}L^{q}_{x'}}\leq C
      h^{-\frac1{p}}\|H_{n}\|_{L^{2}}\lesssim n^{\frac1{p}}\|H_{n}\|_{L^{2}},\
\end{equation}
for all $(p,q)$ satisfying the admissibility condition
$$\frac{2}{p}+\frac{3}{q}=\frac{3}{2}\ ,\ p\geq 2.$$
An identical argument is valid for $\widetilde u$. In fact, for
$\widetilde u$ we shall only need the energy estimate
\begin{equation}\label{ell2}
\| \widetilde u\| _{L^\infty _{x_{1}}L^2_{x'}}\leq C \|\widetilde
H_{l}\| _{L^{2}}\ . \end{equation} Finally, we estimate the
product $u\widetilde u$ as follows. By the Sobolev inequality,
$$\|(1-\Delta)^{-\frac{\alpha}{2}} (u\Tilde u)\|_{L^{2}}
     \leq C\|u\Tilde u\|_{L^{q}},\quad
     \frac 1 q=\frac 1 2+\frac\alpha 4.
    $$
Applying the H\"older inequality we obtain
$$\|(1-\Delta)^{-\frac{\alpha}{2}} (u\Tilde
u)\|_{L^{2}}\leq C\|u\|_{L^{q}_{x_{1}}L^{\frac 4\alpha}_{x'}}
     \|\Tilde u\|_{L^{\infty}_{x_{1}}L^{2}_{x'}}
$$
Noticing that $q<2$ and using the compactness of the support of
$u$, we have
$$\|u\|_{L^{q}_{x_{1}}L^{\frac 4\alpha}_{x'}} \leq C\|u\|_{L^{2}_{x_{1}}L^{\frac 4\alpha}_{x'}}
   .
$$
Applying the Strichartz estimate \eqref{eq.strch} with $p=2$ and
the Sobolev embedding in the $x'$ variables, we obtain
\begin{equation}
\|u\|_{L^{2}_{x_1}(L^{\frac{4}{\alpha}}_{x'})} \leq C
n^{{\frac12-\frac{3\alpha}{4}}} \|u\|_{L^2_{x_{1}}L^{6}_{x'}}\leq
  C n^{{1-\frac{3\alpha}{4}}} \|H_n\|_{L^2}.
\end{equation}
Combining with the $L^\infty L^2$ estimate \eqref{ell2} on $\Tilde
u$, this completes the proof.

\end{proof}

We now come to a quadrilinear estimate on spherical harmonics.

\begin{lemma}\label{quadhs}
Let $\alpha\in ]0,\frac{1}{2}]$ and  $s_0 = 1- \frac{3\alpha}{4}$.
%$\frac{3 \alpha}{4}$,
 There exists $C>0$  such that for any $H_{n_j}^{(j)}, j=1,\cdots ,4$, spherical harmonics on $\mathbb{S}^4$
of degree $n_j$ respectively,
the following quadrilinear estimate holds:
\begin{equation}\label{QQ}
\int _{\mathbb{S}^4}(1-\Delta)^{-{\alpha}} (H_{n_1}^{(1)}
H_{n_2}^{(2)})H_{n_3}^{(3)} H_{n_4}^{(4)} \, dx \leq C
(1+m((n_j))^{s_0}
\prod_{j=1}^4\|H_{n_j}^{(j)}\|_{L^2(\mathbb{S}^4)}.
\end{equation}
\end{lemma}
\begin{proof}
By symmetry, it is sufficient to consider the two cases
$$m(n_1,\cdots ,n_4)= n_1 n_3\quad  ;\quad  m(n_1,\cdots ,n_4)= n_1 n_2\ .$$ In
the first case, the proof follows directly by the Cauchy-Schwarz
inequality and Lemma \ref{BB}. It remains to consider only the
case $m(n_1,\cdots ,n_4)= n_1 n_2$. We use the same idea as in
Lemma \ref{BB} to decompose, if $n_j\geq 1$, each $H_{n_j}^{(j)}$
into  a sum of terms of the form
$$u_j=\chi
_j(x,h_jD_x)H_{n_j}^{(j)}\ ,\ h_j=(n_j(n_j+3))^{-1/2}\ ,\
j=1,2,3,4\ .$$ As before, each $u_j$ may be microlocalized either
into the elliptic zone, in which case we have much stronger
semiclassical estimate \eqref{chi}, in particular an $L^\infty $
bound, or near the characteristic set, and for these terms we can
use the Strichartz type estimate \eqref{eq.strch}. Notice that the
very special case $n_j=0$ can be included into the elliptic case.
Thus we have several possibilities to consider.

If at least two $u_j$'s are microlocalized in the elliptic zone,
then the quadrilinear estimate holds trivially (with $s_0=0$) by a
simple application of the Cauchy-Schwarz inequality.

If $u_3$ or $u_4$ is microlocalized in the elliptic zone, then,
again by the Cauchy- Schwarz inequality,  the quadrilinear
estimate is a consequence of estimate \eqref{eqBB} of Lemma
\ref{BB}, with $\alpha $ replaced by $2\alpha $.

It remains to deal with the cases when only  $u_1$ or $u_2$
 is microlocalized in the elliptic zone, and when all the $u_j$ 's are microlocalized near
 the characterictic set. In both cases, we shall make use of the
 following variant of the Sobolev inequality.

\begin{lemma}\label{sob}
Let $A$ be a pseudodifferential operator of order $-2\alpha $ on
$\mathbb{R}^4$, and let $B$ be a bounded subset of $\mathbb{R}^4$.
For any smooth  function $F$ on $\mathbb{R} ^4$ with support in
$B$, we have the estimate
\begin{equation}
\|A(F)\|_{L^\infty_{x_1}(L^q_{x'})}\leq C
     \|F\|_{L^\infty_{x_1}(L^1_{x'})}
\end{equation}
provided $\frac1q>1-2\alpha/3$.
\end{lemma}
\begin{proof}
The  kernel $K(x,y)$ of $A$ admits an estimate like
\begin{equation}
|K(x,y)| \leq \frac{C}{(|x-y|)^{4-2\alpha}} \leq
\frac{C}{(|x_1-y_1|+ |x'-y'|)^{4-2\alpha}}.
\end{equation}
The claim is then a consequence of Young's inequality in variables
$x'$.
\end{proof}

 By the
self-adjointness of $(1-\Delta)$ the terms to estimate can be
written as follows:
\begin{equation}\label{cof}
I= \left | \int _{\mathbb{S}^4}(u_1 u_2)
 \times (1-\Delta)^{-{\alpha}}(u_3u_4)\, dx\right |\ .
\end{equation}
As in the proof of Lemma \ref{BB} we select a splitting
$x=(x_1,x')$ of the local coordinates such that $u_2,u_3,u_4$ are
solutions of semiclassical evolution equations, and therefore
satisfy Strichartz estimates \eqref{eq.strch}. Using the $L^\infty
$ bound on $u_1$, we have
\begin{equation*}
I \leq C \|H_{n_1}^{(1)}\|_{L^2(\mathbb{S}^4)}\,
\|u_2\|_{L^1_{x_1}( L^{q'}_{x'})} \|(1-\Delta)^{-\alpha}
(u_3u_4)\|_{L^\infty_{x_1}(L^{q}_{x'})},
\end{equation*}
and by Lemma \ref{sob} we obtain
\begin{equation*}
I \leq C \|H_{n_1}^{(1)}\|_{L^2(\mathbb{S}^4)}\,
\|u_2\|_{L^1_{x_1}( L^{q'}_{x'})} \|u_3
u_4\|_{L^\infty_{x_1}(L^{1}_{x'})}
\end{equation*}
provided $\frac{1}{q} > 1- \frac{2\alpha}{3}$.
 H\"older's inequality gives
\begin{equation*}
I \leq C \|H_{n_1}^{(1)}\|_{L^2(\mathbb{S}^4)}\,
\|u_2\|_{L^{2}_{x_1}(L^{q'}_{x'})} \|u_3\|_{L^\infty_{x_1}(
L^{2}_{x'})} \| u_4\|_{L^\infty_{x_1}( L^{2}_{x'})},
\end{equation*}
and, applying estimate \eqref{ell2} on $u_3,u_4$ and estimate
\eqref{eq.strch} with $p=2$ on $u_2$, we obtain
$$I\leq Cn_2^{s}\prod _{j=1}^4\|H_{n_j}^{(j)}\|_{L^2(\mathbb{S}^4)}\,
,$$ with $$s=\max \left(\frac{1}{2},1-\frac{3}{q'}\right )<s_0\
,$$ since $q'$ is arbitrary with $\frac{1}{q'}<\frac{2\alpha}{3}$.

Finally, we treat the  case when all the factors are
microlocalized near the characteristic set. Once again, we select
a splitting $x=(x_1,x')$ of the local coordinates for which
Strichartz estimates \eqref{eq.strch} are valid for each $u_j$. By
H\"older's inequality and Lemma \ref{sob} we have
\begin{equation*}\begin{aligned}
I &\leq C \|u_1u_2\|_{L^1_{x'}( L^{q'}_{x'})}
\|u_3u_4)\|_{L^\infty_{x_1}(L^{1}_{x'})}\\
&\leq C\|u_1\|_{L^2_{x_1}( L^{2q'}_{x'})}
\|u_2\|_{L^2_{x_1}(L^{2q'}_{x'})} \|u_3\|_{L^\infty_{x_1}(
L^{2}_{x'})} \| u_4\|_{L^\infty_{x_1}( L^{2}_{x'})}.
\end{aligned}
\end{equation*}
By estimates \eqref{eq.strch} with $p=2$ on $u_1,u_2$ and
\eqref{ell2} on $u_3,u_4$, we conclude
$$I\leq C(n_1n_2)^s\prod
_{j=1}^4\|H_{n_j}^{(j)}\|_{L^2(\mathbb{S}^4)}\ ,$$ with
$$s=\max\left (\frac{1}{2},1-\frac{3}{2q'}\right )<s_0\ ,$$
since  $q'$ is arbitrary with $\frac{1}{q'}<\frac{2\alpha}{3}$.
This completes the proof.
\end{proof}
\begin{os}\label{Sogge}
It is clear that Lemma \ref{BB} and Lemma \ref{quadhs} extend to
Laplace eigenfunctions on arbitrary compact four-manifolds.
Moreover, a refinement of the study of the elliptic case shows
that, as in \cite{BGT6}, \cite{BGT7}, eigenfunctions can be
replaced by functions belonging to  the range of spectral
projectors of the type ${\bf 1}_{[n,n+1]}(\sqrt {-\Delta}).$
\end{os}
We now come to the main result of this subsection.
\begin{proposizione}\label{Q}
For every $\alpha >0$, for every $s_0>1-\frac{3\alpha }{4}$, the
quadrilinear estimate \eqref{QS} holds on $\mathbb{S}^4$.
\end{proposizione}
\begin{proof}
Let $f_1,\cdots ,f_4$ be functions on $\mathbb{S}^4$ satisfying
the spectral localization property
\begin{equation}\label{Q2}
\mathbf{1}_{\sqrt{1-\Delta} \in [N_j,2N_j]}(f_j) = f_j, \;\;
j=1,2,3,4\ .
\end{equation}
This implies that one can expand
$$f_j=\sum _{n_j}H_{n_j}^{(j)}\ ,$$
where $H_{n_j}^{(j)}$ are spherical harmonics of degree $n_j$, and
where the sum on $n_j$  bears on the domain
\begin{equation}\label{specnj}
N_j/2\leq 1+n_j\leq 2N_j\ . \end{equation} Consequently, the
corresponding solutions of the linear Schr\"odinger equation are
given by
$$u_j(t)=S(t)f_j=\sum _{n_j}{\rm e}^{-itn_j(n_j+3)}H_{n_j}^{(j)}$$
and we have to estimate the expression
\begin{equation*}\begin{aligned}
Q(f_1,\cdots ,f_4,\tau)&=\int_{{\mathbb R}} \int_{\mathbb{S}^4}
\chi (t)\, {\rm e}^{it \tau} (1-\Delta)^{-{\alpha}}
 ( u_1 \overline{u}_2) {u_3} \overline{u}_4 dx dt\\
 &=\sum _{n_1,\cdots ,n_4}\widehat \chi
 (\sum_{j=1}^4{\varepsilon_j}n_j(n_j+3)-\tau)\,
I(H_{n_1}^{(1)},\cdots ,H_{n_4}^{(4)})\ ,
\end{aligned}
\end{equation*}

with $\varepsilon _j=(-1)^{j-1}$ and
$$I(H_{n_1}^{(1)},\cdots
,H_{n_4}^{(4)})=\int_{\mathbb{S}^4}(1-\Delta)^{-{\alpha}}
 ( H_{n_1}^{(1)} \overline{H}_{n_2}^{(2)}) {H_{n_3}^{(3)}}
 \overline{H}_{n_4}^{(4)}\,
 dx\ .
 $$
 Appealing to Lemma \ref{quadhs}, we have, with $s=1-3\alpha /4$,
 $$|I(H_{n_1}^{(1)},\cdots
,H_{n_4}^{(4)})|\leq C\, m(N_1,\cdots ,N_4)^{s}\prod _{j=1}^4\|
H_{n_j}^{(j)}\| _{L^2}\ .$$ Using the fast decay of $\widehat
\chi$ at infinity, we infer
\begin{equation*}\begin{aligned}|Q(f_1,\cdots ,f_4,\tau)|&\leq C\,
m(N_1,\cdots ,N_4)^{s}\sum_{\ell \in \mathbb{Z}}(1+|\ell
|^2)^{-1}\sum _{\Lambda ([\tau]+\ell)}\prod _{j=1}^4\|
H_{n_j}^{(j)}\| _{L^2} \\&\lesssim m(N_1,\cdots ,N_4)^{s}\sup
_{k\in \mathbb{Z}}\sum _{\Lambda  (k)}\prod _{j=1}^4\|
H_{n_j}^{(j)}\| _{L^2}\ ,\end{aligned}
\end{equation*} where $\Lambda  (k)$ denotes the set of $(n_1,\cdots
,n_4)$ satisfying \eqref{specnj} for $j=1,2,3,4$ and
$$ \sum _{j=1}^4\varepsilon _jn_j(n_j+3)=k \ .$$
Now we write
$$\{ 1,2,3,4\} =\{ \alpha ,\beta ,\gamma ,\delta \} $$
with $m(N_1,\cdots ,N_4))=N_\alpha N_\beta$, and we split the sum
on $\Lambda  (k)$ as \begin{equation}\label{split} |Q(f_1,\cdots
,f_4,\tau)|\lesssim m(N_1,\cdots ,N_4)^{s}\sup _{k\in \mathbb{Z}}
\sum _{a\in \mathbb{Z}}S(a)\, S'(k-a) \end{equation} where
\begin{equation*}\begin{aligned}S(a)&=\sum _{ \Gamma (a)}\| H_{n_\alpha }^{(\alpha )}\| _{L^2}\| H_{n_\gamma
}^{(\gamma )}\| _{L^2}\ ;\  S'(a')=\sum _{ \Gamma '(a')}\|
H_{n_\beta
}^{(\beta)}\| _{L^2}\| H_{n_\delta }^{(\delta )}\| _{L^2}\ ,\\
\Gamma (a)&=\{ (n_\alpha ,n_\gamma ): \eqref{specnj}\  {\rm holds\
for}\ j=\alpha ,\gamma ,\ \sum _{j\in \{ \alpha
,\gamma\}}\varepsilon
_jn_j(n_j+3)=a\} ,\\
\Gamma '(a')&=\{ (n_\beta ,n_\delta ): \eqref{specnj}\  {\rm
holds\ for}\ j=\beta ,\delta ,\ \sum _{j\in \{ \beta
,\delta\}}\varepsilon _jn_j(n_j+3)=a'\} .
\end{aligned}\end{equation*}
Now we appeal to the following elementary result of number theory
(see e.g. Lemma 3.2 in \cite{BGT8}).
\begin{lem}\label{count}
Let $\sigma \in \{ \pm 1\}$. For every $\varepsilon >0$, there
exists $C_\varepsilon $ such that, given $M\in \mathbb{Z}$ and a
positive integer $N$,
$$\# \{ (k_1,k_2)\in \mathbb{N}^2: N\leq k_1\leq 2N\ ,\
k_1^2+\sigma k_2^2=M\} \leq C_\varepsilon N^\varepsilon \ .$$
\end{lem}
A simple application of Lemma \ref{count} implies, for every
$\varepsilon >0$,
$$\sup _a\# \Gamma (a)\leq C_\varepsilon N_\alpha ^\varepsilon\ ;\
\sup _{a'}\# \Gamma '(a')\leq C_\varepsilon N_\beta ^\varepsilon\
,$$ and consequently, by a repeated use of the Cauchy-Schwarz
inequality, \begin{equation*}\begin{aligned}\sum _a&S(a)\,
S'(k-a)\leq
C_\varepsilon\, (N_\alpha N_\beta )^\varepsilon \times \\
&\left (\sum _a\sum _{\Gamma (a)}\| H_{n_\alpha }^{(\alpha )}\|
_{L^2}^2\| H_{n_\gamma }^{(\gamma )}\| _{L^2}^2\right )^{1/2}
\left (\sum _a\sum _{\Gamma '(k-a)}\| H_{n_\beta }^{(\beta)}\|
_{L^2}^2\| H_{n_\delta }^{(\delta )}\| _{L^2}^2\right )^{1/2}\\
&\leq C_\varepsilon\, (N_\alpha N_\beta )^\varepsilon \prod
_{j=1}^4\| f_j\| _{L^2}\ ,
\end{aligned}\end{equation*}
where, in the last estimate, we used the orthogonality of the
$H_{n_j}^{(j)}$'s as $n_j$ varies. Coming back to \eqref{split},
this completes the proof. \end{proof}
\begin{os}
Using the remark before the statement of this proposition, the
proof above extends easily to any compact four-dimensional Zoll
manifold (see \cite{BGT8} for more details).
\end{os}

\subsection{Trilinear estimates on the sphere}

\noindent In this subsection, we prove trilinear estimates
(\ref{TS}) on $\mathbb{S}^4$, for every $s_0>1/2$, for zonal
solutions of the Schr\"odinger equation. In view of subsections
2.2 and 2.4, this will complete the proof of  Theorem \ref{t2} and
of Corollary \ref{C1}, by choosing for $G$ the group of rotations
which leave invariant a given pole on $\mathbb{S}^4$.

\vskip 0.25cm \noindent First we recall the definition of  zonal
functions.
\begin{definizione}
Let $d \geq 2$, and let us fix a pole on $\mathbb{S}^d$. We shall
say that a function on $\mathbb{S}^d$ is a zonal function if it
depends only on the geodesic distance to the pole.
\end{definizione}
\noindent  The zonal functions can be expressed in terms of zonal
spherical harmonics which in their turn can be expressed in terms
of classical polynomials (see e.g. \cite{S1}). As in \cite{BGT7},
we can represent the normalized zonal spherical harmonic $Z_{p}$
in the coordinate $\theta$ (the geodesic distance of the point $x$
to our fixed pole) as follows:
\begin{equation}\label{eq.1repr}
    Z_{p}(x)=C(\sin\theta)^{-\frac{d-1}{2}}\left\{
        \cos[(p+\alpha)\theta+\beta]+
        \frac{\mathcal{O}(1)}{p\sin\theta}\right\},\quad
        \frac c p\leq\theta\leq\pi-\frac c p
\end{equation}
with $\alpha,\beta$ independent of $p$, and $C$ uniformly bounded
in $p$. On the other hand, near the concentration points
$\theta=0,\pi$ we can write
\begin{equation}\label{eq.2repr}
    |Z_{p}(x)|\leq C p^{\frac{d-1}{2}},\quad
    \theta\not\in[c/p,\pi-c/p].
\end{equation}
and $\|Z_p\|_{L^2(\mathbb{S}^d)} = 1$.

With this notation, we have the following trilinear eigenfunction
estimates.
\begin{lemma}\label{Z}
There exists a constant $C>0$ such that the following trilinear
estimate holds:
\begin{equation}\label{trilZ}
    \|Z_{p}Z_{q}Z_{l}\|_{L^{1}(\mathbb{S}^{4})}\leq
      C(\min(p,q,l))^{1/2}\ .
\end{equation}
\end{lemma}

\begin{proof}
It is not restrictive to assume that $ p\leq q\leq l$. Moreover,
by Cauchy-Schwarz inequality it is sufficient to  prove
\eqref{trilZ} in the special case $q=l$. Then we have
\begin{equation*}
    \|Z_{p}Z_{q}^{2}\|_{L^{1}(S^{4})}=
    c\int_{0}^{\pi}|Z_{p}(\theta)|Z_{q}(\theta)^{2}(\sin\theta)^{3}d\theta\
    ,
\end{equation*}
where $c$ is some universal constant. We split the interval
$[0,\pi]$ into the intervals $I_{1}=[0,c/q]$, $I_{2}=[c/q,c/p]$,
$I_{3}=[c/p,\pi/2]$ and $I_{4}=[\pi,2,\pi-c/p]$,
$I_{5}=[\pi-c/p,\pi-c/q]$, $I_{6}=[\pi-c/p,\pi]$. Clearly, by
symmetry, it is sufficient to estimate the integral on the first
three intervals $I_{1},I_{2},I_{3}$.

On $I_{1}$ we can use \eqref{eq.2repr} for both harmonics
$Z_{p},Z_{q}$ and the simple estimate $\sin\theta\leq\theta$, and
we obtain
\begin{equation*}
    \int_{0}^{c/q}|Z_{p}|Z_{q}^{2}(\sin\theta)^{3}d\theta\leq
    Cp^{3/2}q^{3}\int_{0}^{c/q}\theta^{3}d\theta\leq
    Cp^{3/2}q^{3}q^{-4}\leq Cp^{1/2}
\end{equation*}
since $q\geq p$.

On the second interval $I_{2}$ we use \eqref{eq.1repr} for $Z_{p}$
and \eqref{eq.2repr} for $Z_{q}$:
\begin{equation*}
    \int_{c/q}^{c/p}|Z_{p}|Z_{q}^{2}(\sin\theta)^{3}d\theta\leq
    Cp^{3/2}\int_{c/q}^{c/p}
    \left(1+\frac 1{q\sin\theta}\right)^{2}d\theta
\end{equation*}
and by the elementary inequality
\begin{equation}\label{eq.elem}
    \left(1+\frac 1{q\sin\theta}\right)^{2}\leq
    C+\frac C{q^{2}\theta^{2}}
\end{equation}
we have immediately
\begin{equation*}
    \int_{c/q}^{c/p}|Z_{p}|Z_{q}^{2}(\sin\theta)^{3}d\theta\leq
    Cp^{3/2}\left(\frac c p-\frac c q+\frac {C}{q^{2}}
    (q/c-p/c)\right)\leq Cp^{1/2}.
\end{equation*}

Finally, in the interval $I_{3}$ we must use \eqref{eq.1repr} for
both harmonics:
\begin{equation*}
     \int_{c/p}^{\pi/2}|Z_{p}|Z_{q}^{2}(\sin\theta)^{3}d\theta\leq
    C\int_{c/p}^{\pi/2}
    \left(1+\frac 1{p\sin\theta}\right)
    \left(1+\frac 1{q\sin\theta}\right)^{2}
    (\sin\theta)^{-3/2}d\theta.
\end{equation*}
Using again \eqref{eq.elem}, the inequality $\sin\theta\geq
C\theta$ on $[0,\pi/2]$, and the fact that $q\geq p$, we have
easily
\begin{equation*}
        \left(1+\frac 1{p\sin\theta}\right)
    \left(1+\frac 1{q\sin\theta}\right)^{2}
    (\sin\theta)^{-3/2}\leq
    C\theta^{-3/2}+Cp^{-3}\theta^{-9/2}.
\end{equation*}
Then integrating on $I_{3}$ we obtain
\begin{equation*}
     \int_{c/p}^{\pi/2}|Z_{p}|Z_{q}^{2}(\sin\theta)^{3}d\theta\leq
     Cp^{1/2}
\end{equation*}
and this concludes the proof.
\end{proof}

We now come to the main result of this subsection, which asserts
that trilinear estimates \eqref{TS} hold for every $s_0>1/2$ on
$M=\mathbb{S}^4$ in the particular case of zonal Cauchy data.

\begin{proposizione}\label{PZ}
Let $s_0 > \frac{1}{2}$ and $\chi \in {\mathcal C}^\infty
_0({\mathbb R})$. There exists $C>0$ such that for any $f_1, f_2,
f_3 \in L^2(\mathbb{S}^4)$ are zonal functions and satisfying
\begin{equation}\label{Pn}
\mathbf{1}_{\sqrt{1-\Delta} \in [N_j,2N_j]}(f_j) = f_j, \;\;
j=1,2,3,
\end{equation}
one has the following trilinear estimate for $u_j(t)=S(t)f_j$,
\begin{equation}\label{TE2}
\begin{aligned}
\sup _{\tau \in \mathbb{R}}\left| \int_{{\mathbb R}}
\int_{\mathbb{S}^4}\right.
 & \left. \chi (t)\,  e^{it \tau}  u_1  u_2 \overline{ u}_3 dx dt  \right| \\
 &\leq C(\min(N_1,N_2, N_3))^{s_0} \|f_1\|_{L^2(\mathbb{S}^4)} \|f_2\|_{L^2(\mathbb{S}^4)}
 \|f_3\|_{L^2(\mathbb{S}^4)}.
\end{aligned}
\end{equation}

\end{proposizione}
\begin{proof} The proof is very similar to the one of Proposition
\ref{Q}. We write
\begin{equation*}
u_j (t)= \sum_{n_j} e^{-itn_j( n_j +3)} c_j(n_j)Z_{n_j},
\end{equation*}
where $n_j$ is subject to the condition \eqref{specnj} and
$$\sum _{n_j}|c_j(n_j)|^2\sim \| f_j\| _{L^2}^2\ .$$
 Thus we can write the integral of the left hand-side of
(\ref{TE2}) as
\begin{equation*}
J=\sum_{n_1,n_2,n_3} \widehat {\chi }(\sum _{j=1}^3\varepsilon
_jn_j(n_j+3) -\tau )\, c_1(n_1)c_2(n_2)\overline {c_3(n_3)}\,
\int_{\mathbb{S}^4}
  Z_{n_1} Z_{n_2} Z_{n_3}dx  \ ,
\end{equation*}
 where $\varepsilon _1=\varepsilon _2=1$ and $\varepsilon _3=-1$. Using the fast
decay of the Fourier transform $\widehat {\chi }$ and the estimate
of Lemma \ref{Z}, we obtain
\begin{equation*}\begin{aligned}
|J|\leq &C\, (\min(N_1,N_2,N_3))^{\frac{1}{2}}\, \sum _{\ell \in
{\mathbb Z}}\frac{1}{1+\ell ^2}\sum _{\Lambda _{[\tau ]+\ell}}
 |c_1(n_1)c_2(n_2)c_3(n_3)|\
  ,\\
  &\lesssim (\min(N_1,N_2,N_3))^{\frac{1}{2}}\, \sup _{k\in \mathbb{Z}}\sum _{\Lambda _k}
  |c_1(n_1)c_2(n_2)c_3(n_3)|\ ,
  \end{aligned}\end{equation*}
where
$$\Lambda _{k}=\{(n_1,n_2,n_3): \eqref{specnj}\ {\rm holds\ for}\ j=1,2,3\ ;\ \sum _{j=1}^3\varepsilon _j
n_j(n_j+3) =k \ \} \ .$$

Suppose for instance that $\min (N_1,N_2,N_3)$ is $N_1$ or $N_2$.
Introducing
$$\Lambda _{k}(n_3)=\{ (n_1,n_2):(n_1,n_2,n_3)\in \Lambda
_{k}\},$$
  we specialize index $n_3$ in the above sum
as \begin{eqnarray*}\begin{aligned} J\leq &C\sup _{k}\sum
_{n_3}|c_3(n_3)|
\left (\sum _{(n_1,n_2)\in \Lambda _k(n_3)}|c_1(n_1)c_2(n_2)|\right )\\
\leq &C\sup _k\left (  \sum _{n_3}|c_3(n_3)|^2   \right
)^{\frac{1}{2}}\left (\sum _{n_3}\bigl (\sum _{(n_1,n_2)\in
\Lambda _k(n_3)}|c_1(n_1)c_2(n_2)|\bigl )^2\right )^{\frac{1}{2}}\\
\leq &C\left ( \sum _{n_3}|c_3(n_3)|^2   \right )^{\frac{1}{2}}
\sup _{k}\left ( \sum _{n_3}[\# \Lambda _{k}(n_3)]\sum
_{(n_1,n_2)\in \Lambda _{k}(n_3)}|c_1(n_1)|^2 |c_2(n_2)|^2\right
)^{\frac{1}{2}}\ .
\end{aligned}
\end{eqnarray*}
To complete the proof, it remains to appeal once again to Lemma
\ref{count}, which yields the estimate
$$\# \Lambda _{\tau ,\ell }(n_3)\leq C_\delta (\min
(N_1,N_2))^\delta \ ,$$ for every $\delta >0$. If $N_3$ is $\min
(N_1,N_2,N_3)$, the proof is similar, by specializing the sum with
respect to $n_1$, say.  \end{proof}

\end{document}